\documentclass[11pt]{amsart}
\usepackage{amssymb,latexsym,amsmath}
\numberwithin{equation}{section}
\newtheorem{definition}[equation]{Definition}
\newtheorem{theorem}[equation]{Theorem}
\newtheorem{proposition}[equation]{Proposition}

\newtheorem{lemma}[equation]{Lemma}

\newtheorem{remark}[equation]{Remark}
 
\def\a{{\alpha}}
\def\b{\beta}

\def\d{\delta}

\def\m{\mu}

\def\w{\omega}

\def\th{\theta}

\def\oa{{\overline \a}}
\def\ob{{\overline \b}}

\def\bC{\mathbb{C}}

\def\bP{\mathbb{P}}
\def\bR{\mathbb{R}}
\def\bZ{\mathbb{Z}}

\def\cF{\mathcal{F}}

\def\fa{\mathfrak{a}}
\def\ff{{\mathfrak{f}}}
\def\fg{{\mathfrak{g}}}
\def\fgl{\mathfrak{gl}}

\def\fp{\mathfrak{p}}

\def\fsl{\mathfrak{sl}_{m+1}}

\def\ft{\mathfrak{t}}

\def\tad{\mathrm{ad}}
\def\td{\mathrm{d}}
\def\tdim{\mathrm{dim}}
\def\tFlag{\mathit{Flag}}
\def\tFub{\mathrm{Fub}}

\def\tGL{\mathrm{GL}}
\def\tker{\mathrm{ker}}
\def\tmod{\textrm{ mod }}
\def\tId{\textrm{Id}}

\def\tneg{\mathrm{neg}}

\def\tred{\mathrm{red}}
\def\tSL{\mathrm{SL}}
\def\tspan{\mathrm{span}}
\def\tSeg{\mathit{Seg}}
\def\l{\ell}
\def\op{\oplus}
\def\ot{\otimes}

\def\ule{\underline{e}}

\def\half{\tfrac{1}{2}}

\def\fourth{\tfrac{1}{4}}

\def\threehalves{\tfrac{3}{2}}

\def\s{\sigma}

\def\xx#1{x^{#1}}

\def\ooo#1#2{\omega^{#1}_{#2}}
\def\oo#1{\omega^{#1}_0}

\def\cf{\mathcal F}

\def\trank{\text{rank}}

\def\BC{\mathbb C}\def\BF{\mathbb F}\def\BO{\mathbb O}
\def\BR{\mathbb R}

\def\BP{\mathbb P}
\def\pp#1{\mathbb P^{#1}}
\def\fa{\mathfrak a}\def\fr{\mathfrak r}
\def\fgl{\mathfrak g\mathfrak l}

\def\fc{\mathfrak c}

\def\pp#1{{\mathbb P}^{#1}}
\def\tdim{\rm dim}
\def\hd{,...,}
\def\ww{\wedge}
\def\upperp{{}^\perp}

\def\inv{{}^{-1}}

\def\cF{{\mathcal F}}\def\cN{{\mathcal N}}

\def\cW{{\mathcal W}}
\def\cO{{\mathcal O}}

\def\BZ{\mathbb Z}

\def\11{\mathbf 1}

\def\fsl{{\mathfrak {sl}}}

\def\ff{{\mathfrak f}}

\def\fg{{\mathfrak g}}
\def\fn{{\mathfrak n}}
\def\fp{{\mathfrak p}}

\def\ft{{\mathfrak t}}
\def\fl{{\mathfrak l}}
\def\l{\lambda}
\def\a{\alpha}

\def\o{\omega}

\def\O{\Omega}
\def\b{\beta}

\def\s{\sigma}

\def\d{\delta}
\def\th{\theta}
\def\m{\mu}
\def\up#1{{}^{({#1})}}

\def\ot{{\mathord{\,\otimes }\,}}
\def\op{{\mathord{\,\oplus }\,}}

\def\ctimes{{\mathord{\times\cdots\times}\;}}

\def\ra{{\mathord{\;\rightarrow\;}}}

\def\La#1{\Lambda^{#1}}
\def\tim{\text{Image}\,}

\def\tann{\text{Ann}\,}
\def\tdim{\text{dim}\,}

\def\tker{\text{ker}\,}
\def\tspan{\text{span}\,}
\def\tmod{\text{ mod }}

\def\trank{\text{rank}\,}

\def\be{\begin{equation}}
\def\ene{\end{equation}}

\def\gp#1{\fg^{\perp}_{#1}}

\def\tii{\widetilde{\mathit{II}}{}^{-2}}

\begin{document}
\title{Fubini-Griffiths-Harris rigidity and Lie algebra cohomology}
\author{J.M. Landsberg \& C. Robles}
\date{\today}
\begin{abstract}
We prove a general rigidity theorem for represented semi-simple Lie groups.  The theorem is used to show that the adjoint variety of a complex simple Lie algebra $\fg$  (the unique minimal $G$ orbit in $\bP\fg$) is extrinsically rigid to third order (with the exception of $\fg=\fa_1$).  

In contrast, we show that the adjoint variety of $SL_3\bC$ and the Segre product $\tSeg(\pp 1\times \pp n)$ are flexible at order two.  In the $SL_3\bC$ example we discuss the relationship between the extrinsic projective geometry  and the intrinsic path geometry.

We extend machinery developed by Hwang and Yamaguchi, Se-ashi, Tanaka and others to reduce the proof of the general theorem to a Lie algebra cohomology calculation.  The proofs of the flexibility statements use exterior differential systems techniques.
\end{abstract}
\maketitle

\renewcommand{\thefootnote}{\arabic{footnote}}

\section{Introduction}
\subsection{History and statement of the problem}

The problem of determining the projective (or extrinsic)
rigidity of varieties $X \subset \bC\bP^N = \bP^N$ dates back     to Monge and has been studied
by Fubini \cite{fub},  
Griffiths and Harris \cite{GH} and others.  The problem may be stated informally as follows: given a  
 homogeneous variety $Z=G/P\subset \BP U=\pp N$ and an unknown
variety $Y\subset \BP W=\pp M$, how many derivatives do we need to take at a general point of $Y$ to determine whether or not $Y$ is projectively equivalent to $Z$? More precisely, there is a sequence of relative differential invariants
of a projective variety $X\subset \BP^N$, defined at a smooth point $x \in X$ (the {\it Fubini forms}, see \S\ref{fubformssect}) that encode the extrinsic geometric information of $X$.
   Let $T_xX$ and $N_xX=T_x\BP^N/T_xX$ denote the tangent and normal spaces to $X$ at $x \in X$.  Then the $k$-th Fubini form $F_{k,x}$ at $x$ is an element of $S^kT^*_xX\ot N_xX$ modulo an equivalence relation.
If we (locally) express $X$ as graph in $\bP^N$ over its embedded tangent space, then the Fubini forms $F_{k,x}$ are defined by the degree $k$ coefficients in the Taylor series expansion of $X$ about $x$, modulo the equivalence relation given by other choices of first-order adapted local coordinates.

\begin{definition} 
Let $X\subset \BP U$ be a   projective variety and
let $x\in X$ be a general point.  
\newcounter{cl}
\begin{list}{$\circ$}
  {\usecounter{cl}
   \setlength{\leftmargin}{15pt}
   \setlength{\labelwidth}{10pt}
   }  
\item 
If   $Y \subset \bP W$ is another variety
such that for a general point $y\in Y$, there exist    bijective linear maps $\phi_T: T^*_yY\ra T^*_xX$, $\phi_N: N_yY\ra N_xX$,   such that the induced linear maps
$$  (\phi_T)^{\circ \ell}\ot \phi_N: S^{\ell}T^*_yY\ot N_yY
\ra S^{\ell}T^*_xX\ot N_xX$$ take  (equivalence classes of) Fubini forms of $Y$ to (equivalence classes of) Fubini forms of $X$ for all $\ell\leq k$, then we say that \emph{ $Y$ agrees with $X$ to order $k$}.  
\item 
We say \emph{$X$ is rigid at order $k$} if whenever
a variety  $Y$ agrees with $X$ to order $k$, there exists a linear map $\Phi: W\ra U$ such that $\Phi(Y)=X$.
\item 
When $X$ is not rigid at order $k$, but the set of distinct $Y$ agreeing with $X$ to order $k$ is finite dimensional, we say \emph{$X$ is quasi-rigid at order $k$}.
\item
Otherwise we say \emph{$X$ is flexible at order $k$}.
\end{list}
\end{definition}

\noindent{\it Remark.} In this paper we discuss
the rigidity problem in $\bP U$.  It is also natural
to study rigidity when the ambient space is a homogeneous
variety $G'/P'$ and the model variety is a $G$-variety.  See, for example, the study of Schubert varieties in compact Hermitian symmetric spaces in \cite{BrSchubert}.  

None of the examples in this paper are quasi-rigid.  The paper \cite{BrSchubert} contains quasi-rigid examples for a closely related problem.
\smallskip

Projective space embedded as a linear subspace $\bP^n \subset \bP^N$ is easily
seen to be rigid at order two. Fubini \cite{fub} proved that the smooth
$n$-dimensional quadric hypersurface is rigid at order three  when $n>1$.  When $n=1$ Monge proved that the quadric is rigid at order five (cf. \cite[Ex. 2.26]{Lci}). Griffiths and Harris \cite{GH} conjectured
that the Segre variety $\tSeg(\pp 2\times \pp 2)\subset\pp 8$
is rigid at order two, and the conjecture was proven  in \cite{Lrigid}. Then in
\cite{Lchss} it was shown
that any irreducible, rank two, compact Hermitian symmetric space (CHSS)
in its minimal homogeneous embedding (except for the quadric hypersurface) 
is rigid at order two.   Hwang and
Yamaguchi \cite{HY} then solved the rigidity problem for all homogeneously embedded irreducible CHSS: other than the quadric
hypersurface,
an irreducible, homogeneously embedded CHSS with
osculating sequence of length $f$ is rigid at
order $f$. (The length of the osculating filtration of a CHSS in its minimal homogeneous embedding is equal to its rank.)
 
The next simplest homogeneous varieties are the {\it adjoint varieties} 
$Z_{\tad}^G\subset \BP\fg$, the {\it homogeneous complex contact manifolds}.  These are the unique closed orbits in $\bP\fg$ of the adjoint action of the associated complex simple Lie group $G$.  When $G=C_n$,   $Z_{\tad}^{C_n}=v_2(\pp{2n-1})$ is the quadratic Veronese variety which is CHSS for the larger group $A_{2n-1}$. In \cite{Lrigid} $v_2(\pp{m})$, $m>1$, was shown to be rigid at order three.  This variety is flexible at order two since its second fundamental form is generic.
The quadratic Veronese varieties have vanishing third and fourth order Fubini forms. In \cite{LMmagic} it was observed that the third and fourth order Fubini invariants are non-zero for the other adjoint varieties.  This led the authors to speculate that the other adjoint varieties would be rigid at order four, but not three.    Robles \cite{robles} then showed that $Z_\tad^{A_n}$ is rigid at order three.  (Note that, with the exceptions of $\fa_n$ and $\fc_n$, the adjoint representations of the complex simple Lie algebras are fundamental.)  

\subsection{Results}
\begin{theorem}\label{adjointthm}
The adjoint variety of a complex simple Lie algebra $\fg$, $\fg \not= \fa_1$, is rigid at order three.
\end{theorem}
\noindent
Theorem \ref{adjointthm} is a consequence of  Theorem \ref{bigthm} ($\fg \not= \fa_n$),  Lemma \ref{adjointlemma} 
(when $\fg = \fa_n$, $n>2$) and Proposition
\ref{prop:Fub3=>-1}. The case of $A_2$ was resolved in
\cite{robles} and we give a different proof in \S\ref{sl3end}.

\begin{theorem}\label{thm:I0J0}
Assume $n,r>1$.
\noindent
\begin{list}{\emph{(\alph{cl})}}
  {\usecounter{cl}
   \setlength{\leftmargin}{25pt}\setlength{\itemsep}{5pt}
   \setlength{\labelwidth}{20pt}
   } 
\item 
The Veronese varieties $v_d(\pp n)$ are rigid at order $d+1$.
\item
The Veronese embeddings of the quadric hypersurfaces $v_d(Q^n)$ are rigid at order $2d+1$.
\item 
For integers $1 = a_1 \le a_2 \le \cdots \le a_r$, the Segre variety $\tSeg( \bP^{a_1} \times \bP^{a_2} \times \cdots \times \bP^{a_r} )$ is rigid at order $r+1$.
\end{list}
\end{theorem}
\noindent
Theorem \ref{thm:I0J0} is proven in \S\ref{resta}.  (The proof references calculations in \S\ref{compsect}.)

\medskip

\noindent{\it Remark.}  If {\it all} the integers in Theorem \ref{thm:I0J0} (c) satisfy $a_i>1$,  then $\tSeg$ is rigid at order $r$, c.f. \cite{HY}.  The third-order rigidity of $\tSeg(\bP^1 \times \bP^1)\subset\bP^3$  is a consequence of Fubini's theorem \cite{fub},
as it is a quadric hypersurface.
\smallskip

We apply standard exterior differential
systems (EDS) techniques, guided by representation theory, to establish the following two flexibility results in \S\ref{p1pnsect} and \S\ref{sl3sect}, respectively.   We remark that the calculations for Theorem \ref{thm:flexsl3} are quite long.


\begin{theorem}\label{thm:flexseg}The Segre variety $\tSeg(\pp 1\times \pp n)\subset \BP^{2n+1}$,  $n>1$, is flexible at order two.
It is rigid at order three.
\end{theorem}

Theorem \ref{thm:flexseg} is proved in \S\ref{p1pnsect}.  In \cite{Lchss} it was mistakenly remarked that
$\tSeg(\pp 1\times \pp n)$ was rigid at order two. The
reasons for this error are explained in \S\ref{p1pnsect}.
Notice that Theorem \ref{thm:flexseg} implies that Theorem \ref{thm:I0J0} (c) is sharp, and the rigidity statement is a consequence of Theorem \ref{thm:I0J0} (c).  

\begin{theorem}\label{thm:flexsl3}
The adjoint variety $Z_\tad^{A_2}=\tFlag_{1,2}(\BC^3)$
is flexible at order two. The general integral manifold
gives rise to a non-flat path geometry. There
are also integral manifolds that are projectively
inequivalent to the homogeneous model and
give rise to flat path geometries. 
\end{theorem}

Theorem \ref{thm:flexsl3} is proved in \S\ref{sl3sect}.
We actually prove a stronger flexibility result. In \S\ref{rootgradingsect} we  define the $(I_{-1},J_{-1})$ system associated to any homogeneously embedded homogeneous variety. (This system is more restrictive than the second order Fubini system.)
A parabolic $P\subset G$ determines a $\BZ$-grading of $\fg$
and all $\fg$-modules, including $\fgl(U)$, where
$G/P\subset \BP U$. For $p=0,-1$, we define the $(I_p,J_p)$
exterior differential system  by restricting the the component of the Maurer-Cartan
form of $GL(U)$ taking values in the
$q$-th homogeneous part of $\fgl(U)$ to take values in $\fg \subset \fgl(U)$, for $q\leq p$.  See \S\ref{rootgradingsect} for a 
detailed description.

While the $(I_p,J_p)$ systems are natural from the point of view of
representation theory, except in the case of generalized cominuscule representations,
they do not lead one to Lie algebra cohomology. To rectify this, we define
a class of {\it filtered EDS} as follows:
Define $(I^\textsf{f}_{p},\Omega)$ to be the {\it $(p+1)$-filtered EDS} on $GL(U)$  with the
same independence condition as the $(I_p,J_p)$ system  but only requiring integral
manifolds $i: M\ra GL(U)$ to satisfy $i^*(\o_{\fg^{\perp}_{p-s}}|_{T^{-1-s}M})=0$ where
the filtration $T^{-1}M\subset T^{-2}M\subset \cdots \subset T^{-k}M=\cdots = TM$ on $TM$ is induced from the grading on $TGL(U)$. See \S\ref{filtersect} for more details.

For the adjoint varieties, the $(I_{-1},J_{-1})$ system implies second order agreement plus partial third order agreement.  The integral manifolds
are bundles $\cF \subset GL(U)$ over varieties $Y \subset \bP U$.  The EDS imposes the condition that the base varieties $Y$ all have contact hyperplane distributions, and that they   support an intrinsic parabolic geometry.  When $G={SL_3}\bC$ the parabolic   geometry is the incidence space for a {\it path geometry} in the plane.

The variety $Z^{SL_3}_\tad = \tFlag_{1,2}(\BC^3)$ is the homogeneous
model for the incidence space of a path geometry  in the plane. This
geometry has been studied extensively.  The equivalence problem was solved by Cartan \cite{cartan}.
Intrinsically, on a three-fold $X$ equipped
with two foliations by curves
whose tangent lines span a contact distribution, there are two differential invariants,
$J_1$ and $J_2$, which measure the failure of $X$ to
be locally equivalent to the homogeneous model.
We show that there are three distinct types
of integral manifolds.   The first class consists only
of the standard homogeneous model.
The second class depends (in the language of Cartan) on
four functions of one variable. Despite this
flexibility, {\it all the examples in the
second class are intrinsically flat}. 
(One can think by analogy of the surfaces of
zero Gauss curvature in Euclidean 3-space, which depend on functions of one variable.)
The third class of integral manifolds depends roughly on
two functions of two variables. Here the
intrinsic invariants are nonzero.

As mentioned above, the  $(I_{-1},J_{-1})$ system
is more restrictive than specifying
second order agreement.  In the case of the adjoint varieties
it is less restrictive than specifying third
order agreement.
Theorems \ref{adjointthm} and \ref{thm:I0J0} are corollaries
of Theorem \ref{bigthm} below.  Loosely speaking, adjoint
varieties for simple Lie algebras other than
$\fsl_2$ and $\fsl_3$ are rigid to an order between two
and three. In \cite{LMmagic}
it was observed
that adjoint varieties are  not   rigid to order two,  
but    only one non-isomorphic variety that
agreed to order two was known: a \lq\lq parabola\rq\rq\
which had the same second fundamental form
of the corresponding adjoint variety and
all other differential invariants zero. 
These parabolas are not equipped with a contact
structure.

\begin{definition}
A homogeneous variety $Z=G/P\subset\BP U$ is \emph{rigid for the $(I_{p},J_p)$ system}
(or \emph{rigid for the $(I_{p}^f,\Omega)$ system} that is defined in \S\ref{filtersect})  if the only integral manifolds of the system on $SL(U)$
correspond to conjugates of $G\subset SL(U)$.  
\end{definition}
\noindent
In this case the underlying projective varieties are all 
projectively equivalent 
to $Z$. Any variety rigid for the
$(I_{p},J_{p})$ system (resp. the $(I_p^\textsf{f},\Omega)$ system) is automatically
rigid for the $(I_{p+1},J_{p+1})$ system (resp. the $(I_{p+1}^f,\Omega)$ system).
Moreover, any variety rigid for the $(I_p^\textsf{f},\Omega)$ system  is automatically
rigid for the $(I_{p },J_{p })$ system.

Let $G$ have rank $r$, let $I\subset \{ 1\hd r\}$ and
write $P=P_I\subset G$ for the parabolic subgroup obtained by deleting negative root spaces corresponding to roots having a nonzero coefficient on any of the simple roots $\a_{i}$, $i \in I$.

\begin{theorem}\label{bigthm}
Let $G$ be a complex semi-simple group and let $Z=G/P\subset \BP U$ be a homogeneously embedded homogeneous variety
(i.e., the orbit of a highest weight line in $\BP U$). 
Assume that no factor of $Z$ corresponds to a quadric
hypersurface or $A_n/P_I$, with $1$ or $n$ in $I$.  Then $Z$ is rigid for the $(I_{-1}^f,\Omega)$-system.
\end{theorem}
\noindent
The hypotheses above exclude the adjoint variety $Z^{A_n}_\tad$, where $I = \{1,n\}$.  Nonetheless,
\begin{theorem}\label{anadjthm}  The adjoint variety
$Z_{\tad}^{A_n}\subset \BP (\fa_n)$, for $n>2$, is  rigid for the $(I_{-1}^f,\Omega)$-system.
\end{theorem}

\begin{theorem}\label{vdpna2thm}  The Veronese variety
$v_d(\pp n)\subset \BP (S^d\BC^{n+1})=\bP U^{A_n}_{d\o_1}$,
$n\geq 2$, 
and the adjoint variety $Z_{\tad}^{A_2}\subset \BP \fa_2$  
are rigid for the  $(I_0,J_0)$ system.
\end{theorem}
\noindent
Theorems \ref{bigthm}, \ref{anadjthm} and \ref{vdpna2thm}
are proven in \S\ref{compsect}.

\subsection{Reduction to Lie algebra cohomology}\label{proofmeth}
Given a homogeneous variety $Z=G/P\subset \BP U=\pp N$ and
a point $x\in Z$, one obtains two filtrations
of $U$. The first is  the {\it osculating filtration} (see \S\ref{oscfiltsect})
which any manifold in $\bP^N$ has.  This filtration corresponds
to the spans of successive derivatives at $x$ of curves on
$Z$. The second is the {\it grading filtration} which
is given in terms of Lie algebra data (see \S\ref{rootgradingsect}).
The two filtrations coincide if and only if $X$ is a homogeneously
embedded compact Hermitian symmetric space (CHSS).
There are natural exterior differential systems (EDS)
on $GL(U)$  associated
to each of these filtrations. 

In \cite{HY}, Hwang and Yamaguchi used
methods of   Tanaka   and Se-ashi \cite{seashi}
as developed in \cite{SYY} to establish an extrinsic rigidity result
for CHSS subject to partial vanishing of the first Lie algebra
cohomology groups. The methods of \cite{HY} were translated to
the language of EDS in \cite{Lima}. In brief, for CHSS, if one quotients
the kernel of the Spencer differential by admissible normalizations in the fiber one arrives naturally at the Lie algebra cohomology group $H^1(\fg_{-},\fg\upperp)$.

Two problems present themselves in generalizing
beyond the case of CHSS. First, the Fubini systems
and the EDS induced by the Lie algebra grading no
longer coincide. Second, even if one ignores the
Fubini system and restricts attention to the EDS
induced by the Lie algebra gradings, {\it the 
  Spencer differential no longer coincides
with the Lie algebra cohomology differential.}
We resolve the first problem  
for adjoint varieties  by proving that
integral manifolds of the third order
Fubini system must also be integral manifolds of
the $(I_{-1},J_{-1})$ system. 
In a planned sequel, we hope to resolve this problem
for general $G/P$. (Resolving the problem
in individual cases is straightforward, but
tedious using our current technology. For 
example, we explain how to show   that the
$78$-dimensional  variety
$E_8/P_1\subset\pp{3874}$ is rigid to order five in \S\ref{e8sect}.) 
The second difficulty is resolved 
in complete generality in this paper by introducing the filtered EDS
$(I_p^\textsf{f},\Omega)$ and simultaneously examining the torsion equations and
prolongation in lowest nontrivial homogeneous degree (see our main Theorem \ref{mainthm}):

\smallskip

\begin{quote}
{\it The Spencer differential   of the filtered system combined with the torsion
equations  in their first
nontrivial homogenous degree $d$ is the degree $d$ homogeneous component of the Lie algebra cohomology differential $\partial^1$.}
\end{quote}
\smallskip

\noindent We hope this result will be useful in the study
of other EDS with symmetry.

\medskip

\noindent{\it Remark.} The second order Fubini system induces the structure of a {\it parabolic geometry} on integral manifolds.  When the grading on the tangent space is two-step, the third order Fubini system forces the parabolic geometry to be regular (e.g. one has a multi-contact structure on $T_{-1}$).  The differential invariants that measure the failure
of a parabolic geometry to be locally flat take
values in  the Lie algebra cohomology groups $H^2(\fg_{-},\fg)$  (see, e.g., \cite[\S 4]{cappara}).    It may be interesting to compare these intrinsic and extrinsic geometries and their invariants in greater detail.

\subsection{Overview}
 
In \S\ref{edssect} we review basic definitions from exterior differential systems (EDS) and introduce the notion  of {\it reduced  prolongations}.  In \S\ref{fubinisect} we review the Fubini forms and introduce natural EDS (the {\it Fubini systems}) for studying rigidity of projective varieties. 
We then introduce the 
$(I_0,J_0)$ and $(I_{-1},J_{-1})$ EDS  in \S\ref{rootgradingsect}. 
 In \S \ref{motexsecta} we show the $(I_{-1},J_{-1})$ system for
adjoint varieties is implied by the third order
Fubini system. 
 
 We define the filtered EDS $(I_p^\textsf{f},\Omega)$ in \S\ref{filtersect}. We then
review Lie algebra cohomology and Kostant's theory
in \S\ref{liealgcohsect} and \S\ref{kossect}, respectively.  
In \S\ref{filtersect} we reduce proving rigidity for the
$(I_p^\textsf{f},\Omega)$ system to proving that certain Lie algebra cohomology groups are zero.  In \S\ref{compsect} we identify  groups and modules for which these cohomology groups vanish.

We establish flexibility
of $\tSeg(\pp 1\times \pp n)$ and $Z^{A_2}_{\tad}$ respectively
 in \S\ref{p1pnsect} and \S\ref{sl3sect}.  In \S\ref{sec:contact} we   discuss the intrinsic geometry of the integral manifolds modeled on $Z^{A_2}_{\tad}$.
Finally, in \S\ref{restsect} we complete the
proofs of the remaining rigidity results.

\subsection{Generalizations} 
As discussed above, one could  
prove further  rigidity results by
having a better understanding of the Fubini forms.
More precisely, one would need to determine to what
order the Fubini form of a homogeneous variety must be specified in to fix the negative part of the Maurer-Cartan
form in the Lie algebra grading. We plan
to address this question in a subsequent paper.

\smallskip

All our systems have immediate generalizations to systems for
homogeneous submanifolds $G/P\subset G'/P'$ of homogeneous
varieties.  (That is, $\bP U$ may be replaced by the more general 
$G'/P'$.)  Hong \cite{hong1,hong2} has results in this setting for three-step gradings.  The machinery we have developed here will permit a
more general study.

For the study of general Schubert varieties (see, for example, \cite{BrSchubert}), or more generally $G$-subvarieties of a homogeneous variety,  there is another obstruction to applying 
Kostant's theory.  Kostant's method of calculating
$H^1(\fn,\Gamma)$ only applies if $\Gamma$ is a $\fg$-module 
and $\fn\subset \fp\subset \fg$.  It may be possible to extend Kostant's theory to apply to systems
 for varieties modeled on Schubert varieties by finding a sufficiently invariant inner product on the relevant $\fn$ and $\Gamma$
to enable a \lq\lq Hodge type\rq\rq\ theorem.

\subsection{Notational conventions}
We work over the complex numbers throughout, all functions are holomorphic functions,
manifolds complex manifolds etc... (although much of the theory carries over to $\BR$,
with some rigidity results even carrying over to the $C^{\infty}$ setting). In particular
the notion of a {\it general point} of an analytic manifold makes sense, which is a point
off of a finite union of analytic subvarieties.
We use the labeling and ordering of roots and
weights as in \cite{bou}. For subsets $X\subset \BP V$, $\hat X\subset V$ denotes the corresponding cone.  For a manifold
$X$, $T_xX$ denotes its tangent space at $x$. For a submanifold $X\subset \BP V$, $\hat T_xX = T_p\hat X\subset V$, 
denotes its affine tangent
space, and $p \in \hat x =: L_x$. In particular, $T_xX=\hat x^*\otimes \hat T_xX/\hat x$. We use the summation
convention throughout: indices occurring
up and down are to be summed over. If $G$ is semi-simple of rank $r$, we write
  $P=P_I\subset G$ for the parabolic subgroup obtained by deleting negative root spaces corresponding to roots having a nonzero coefficient on any of the simple roots $\a_{i_s}$, $i_s \in I\subset\{1\hd r\}$.

\subsection{Acknowledgments} 
We thank
R. Bryant, B. Doubrov, J. Hwang,  S. Kumar, 
and K. Yamaguchi for useful discussions. The essential step for the transition
to Lie algebra cohomology was resolved while A. Cap was
visiting us.  Remarkably Cap needed almost
identical machinery for his work (see \cite{cap}) and, after a long
day of discussions, both problems were resolved independently.
This paper has benefitted tremendously from our
conversations with him. In particular, the reformulation of our results in
terms of filtered EDS was suggested by Cap. 

\section{Linear Pfaffian systems (EDS)}\label{edssect}

A {\it linear Pfaffian
exterior differential system} (EDS) on a manifold $\Sigma$
is a flag of  subbundles of the cotangent bundle to $\Sigma$,
$I\subset J\subset T^*\Sigma$ such that the map
$I\ra \La 2(T^*\Sigma/J)$ given by $\th\mapsto d\th\tmod J$ is zero. An {\it integral manifold}
is an immersed submanifold $i: N\ra \Sigma$,
where $n=\tdim N=\trank(J/I)$, such that
$i^*(I)=0$ and $i^*(\La n(J/I))$ is non-vanishing.

A standard example is 
the space of one-jets $\Sigma=J^1(\BR,\BR)$ with
coordinates $(x,y,p)$ and $I=\{ dy-pdx\}$
$J=\{ dy-pdx,dx\}$. Integral manifolds are the 
lifts $\{(x,f(x),f'(x))\}$ of graphs of arbitrary differentiable functions $f:\bR\to\bR$.

Fixing a general point $x\in \Sigma$ we set
$V=(J/I)_x^*$, $W=I_x^*$.
Calculating $dI\tmod I$ one defines
the {\it tableau} $A\subset W\ot V^*$
which (in the absence of torsion)
essentially  consists
of the space    of $n$-planes in $T_x\Sigma$
on which the   forms
in $I$ vanish; equivalently, the possible first
order Taylor series of any integral manifold (expressed as a graph) through
$x$.
Because the systems in this paper are modeled on homogeneous spaces,  $V$, $W$ and $A$ will be independent of the base point $x\in \Sigma$: we suppress reference to it.

Calculating $dI \tmod I$ also yields the {\it torsion}, which is
an element of $W\ot \La 2V^*/\tilde\d(A\ot V^*)$, where $\tilde\d$ is the 
Spencer differential defined below.
If there are no $n$-planes in $T_x\Sigma$
on which the forms in $I$ and their derivatives vanish at $x$, one
says the system has {\it torsion} at $x$. In this
case one must restrict to   the submanifold
$\Sigma'\subset \Sigma$ where there is no torsion.
 
Consider $f\in W\ot V^*$ as a linear map $f: V\ra W$ and define 
$$\tilde\d: (W\ot V^*)\ot V^*\ra W \ot \La 2 V^*
$$
by
   $\tilde\d (f\ot \xi)=df\ww\xi$, and
$A^{(1)}:=\tker \tilde\d|_{A\ot V^*}$,
  the {\it prolongation of $A$}.
Cartan developed a test to determine the \lq\lq
size\rq\rq\ of the space of local integral manifolds
by comparing a crude estimate
for the dimension of the
space of  admissible second order terms in the Taylor
series (which is $A\up 1$) with the actual dimension.
When this test fails, one {\it prolongs} the system, defining
a new system on $\Sigma\times A\up 1$.

The systems we will be dealing with   all have
at least one integral manifold, the corresponding homogeneous model.  These systems may be grouped into three distinct categories.
\begin{definition}
We will say a system is \emph{rigid} if all integral
manifolds are equivalent to the homogeneous model.
We will say a system is \emph{quasi-rigid} if
some prolongation of its tableau is zero and
the system is not rigid. (In other words, some
prolongation of the system is a Frobenius system.) In this
case the space of local integral manifolds is finite
dimensional. Otherwise we will say the system
is \emph{flexible}. 
\end{definition}
\noindent
In the case that the system is flexible, the space of
local integral manifolds is infinite dimensional -- the integral manifolds will be determined by some space of free functions, as in the  
example of $J^1(\BR,\BR)$ above.

The {\it first derived system} $I'\subset I$ is the
maximal subbundle such that $dI'\equiv 0\tmod I$
holds. We may and will suppress it in calculations.
(See \cite[p.216]{EDSBK} or \cite[p.45]{BCG3}.)

A vector field $X\in \Gamma(T\Sigma)$ is a
{\it Cauchy characteristic vector field} 
for the linear Pfaffian system $(I,J)$ if $\a(X)=0$ for
all $\a\in I$. Cauchy characteristic vector fields correspond
to (infinitesimal) symmetries of the EDS. By
adding differential forms dual to the Cauchy characteristics to $J$, 
one obtains a new system whose integral manifolds are in natural 
correspondence with those of the original system.
(See \cite[p.205]{EDSBK}, \cite[p.36]{BCG3}.)

\subsection{Reduced prolongations}\label{redprolongsect}
In our case   $\Sigma=GL(U)$, the group of invertible
linear maps $U\ra U$,
and $I$ will consist of a component of the Maurer-Cartan
form of $GL(U)$. We are interested in submanifolds
of $\BP U$; the integral manifolds of our EDS will be frame
bundles $\cF \subset GL(U)$ over these submanifolds. Consider the subgroup of $\cN' \subset GL(U)$ that 
preserves the EDS: it will carry integral manifolds to integral manifolds.
Let $\cN\subset \cN'$ denote the subgroup that preserves the base
submanifold in $\BP U$.  

The prolongation coefficients in $A^{(1)} \subset A \ot V^*$ vary under the action of $\cN$.  These variations are the {\it admissible normalizations} of the prolongation coefficients.  The admissible normalizations may be realized as the image of a map $\delta: \fn \to A \ot V^*$, where $\fn$ is the Lie algebra of $\cN$.  We exploit this group action by defining a
{\it reduced prolongation} of our EDS: 
$$
A_{\tred}^{(1)} =A\up 1/\tim (\d: \fn\ra A\ot V^*) \, .
$$

As far as determining integral manifolds 
via the Cartan algorithm is  concerned,
normalizing by the action of $\cN$ is of
limited practical effect, but {\it for the EDS discussed in this paper,
the normalization corresponds to  the denominator of a Lie algebra cohomology group}.

\section{The Fubini systems}\label{fubinisect}
\subsection{The Fubini forms}\label{fubformssect}
Let 
$X^n\subset \BP U=\pp N$ be
a submanifold and let $x\in X$. Take linear local
coordinates on $\BP U$ about $x$ adapted so that
$X$ is locally graphed over its (embedded) tangent plane at $x$.
Expanding the functions   in Taylor series, 
and collecting the homogeneous terms in each degree,
one obtains a series of tensors 
${\underline F}_k={\underline F}_{k,X,x}\in (\hat x^{\ot^{k-1}}) \ot S^kT^*X
\ot N_xX$.  Here $N_xX=T_x\BP U/T_xX$ denotes
the normal space. If the local
linear coordinates are
$(\xx 1\hd \xx n,\xx{n+1}\hd \xx N)=(\xx\a,\xx\mu)$
so that $X$ is locally given by equations
$\xx\mu=f^{\mu}(\xx\a)$ then
$$
{\underline  F }_{k,X,x}= (-1)^k \, 
\frac{\partial^kf^{\mu}}{\partial \xx{i_1}  \cdots 
\partial\xx{i_k}} \,  dx^{i_1}\circ \cdots\circ  d\xx{i_k}
\ot \frac\partial{\partial\xx\mu}
$$
The ${\underline F_k}$  depend
on the choice of linear coordinates. To rectify this
ambiguity, let $\pi: \cF^1\ra X$ denote the bundle
of first order adapted frames, that is
$$
\cF^1=\big\{ (e_0\hd e_N)\in GL(U) \mid [e_0]\in X, \
\hat T_{[e_0]}X=\tspan\{e_0\hd e_n\} \big\} \subset GL(U) \, .
$$
Here the $e_A$ are to be considered as column vectors.
The fiber over $x$ is isomorphic to the space of
first order adapted linear coordinates based at $x$.
Thus one obtains well defined
$F_k\in \Gamma(\cf^1,\pi^*(L^{k-1} \ot S^kT^*X\ot NX))$
which are called the {\it Fubini forms} \cite{EDSBK}. Here $L$ is the line bundle
$\cO_{\BP U}(-1)$, the
bundle with fiber $L_x=\hat x$.  We use
the notation $L^k := L^{\ot k}$.  The Fubini forms determine $X$ up
to projective equivalence.

\subsection{The Maurer-Cartan form}
\label{sec:mc}
Let $G$ be a  Lie group, let 
$\o\in \O^1(G,\fg)$ denote the {\it Maurer-Cartan form}
of $G$.  By definition $\o$ is the unique left-invariant
differential form such that $\o_{\tId}:T_{\tId}G\simeq \fg\ra\fg$
is the identity map.   It satisfies the {\it Maurer-Cartan equation} 
\begin{equation}\label{eqn:mc}
  d\o \ = \ -\left[\o,\o\right] \, .
\end{equation} 
To be precise, given two $\fg$-valued 1-forms $\rho$ and $\eta$, and two tangent vectors $u,v$, we define 
$$ \left[ \rho , \eta \right](u,v) \ := \ 
   \left[ \rho(u) , \eta(v) \right] - \left[ \rho(v) , \eta(u) \right] \, .
$$
In practice, $\rho$ and $\eta$ will be components (of a fixed homogeneous degree) of the Maurer-Cartan form.

If $G\subseteq GL(U)$ is
a matrix Lie group,  the Maurer-Cartan equation
may be written $d\o=-\o\ww\o$: if $\o=(\ooo AB)$, with $\ooo AB$  the entries of
the matrix, we have
$(d\o)^A_B=-\ooo AC\ww\ooo CB$. An essential
feature of the Maurer-Cartan equation is that differentiation is reduced to an
algebraic operation.

\subsection{Second order Fubini systems}\label{sec:Fub2}
Fix vector spaces $L,T,N$ of dimensions $1,n,a$ and
fix an element $F_2\in S^2T^*\ot N\ot L$. Let $U = L\op T\op N$, and let $\o\in \O^1(GL(U),\fgl(U))$ denote the Maurer-Cartan form. 
Write $\fgl(U)= (L \op T\op N)^*\ot (L\op T\op N)$
and let, for example,  $\o_{L^*\ot T}$ denote the component
of $\o$ taking values in $L^*\ot T$.
We have
\be
\o=\begin{pmatrix}
\o_{L^*\ot L}&\o_{T^*\ot L}&\o_{N^*\ot L}\\
\o_{L^*\ot T}&\o_{T^*\ot T}&\o_{N^*\ot T}\\
\o_{L^*\ot N}&\o_{T^*\ot N}&\o_{N^*\ot N}
\end{pmatrix} \, .
\ene
Writing the Maurer-Cartan equation
component-wise yields, for example,
$$
d\o_{L^*\ot T}=-
\o_{L^*\ot T}\ww\o_{L^*\ot L}
-\o_{T^*\ot T}\ww\o_{L^*\ot T}
-\o_{N^*\ot T}\ww\o_{L^*\ot N}.
$$
 
Given $F_2 \in L \ot S^2T^* \ot N$, 
the {\it second order Fubini system} for $F_2$ is
$$I_{\tFub_2}=\{\o_{L^*\ot N}, \ \o_{T^*\ot N}-F_2(\o_{L^*\ot T})\},\ \ 
 J_{\tFub_2}=\{I_{\tFub_2}, \ \o_{L^*\ot T}\}.
$$
Its integral manifolds are submanifolds $\cF^2 \subset GL(U)$ that
are adapted frame bundles of submanifolds $X\subset \BP U$
having the property that at a point $x\in X$, the projective
second fundamental form $F_{2,X,x}$ is equivalent to $F_2$.
(The tautological system for frame bundles $\cF^1$ of
arbitrary $n$ dimensional submanifolds
is given by $I=\{\o_{L^*\ot N}   \}$, 
$ J=\{I, \o_{L^*\ot T}\}$.)

Let $R\subset GL(L)\times GL(T)\times GL(N)$
denote the subgroup stabilizing $F_2$ and
let 
$$\fr\subset (L^*\ot L)\op (T^*\ot T)\op 
(N^*\ot N) =:\fgl(U)_{0,*}
$$
denote its subalgebra.  These are the elements of $\fgl(U)_{0,*}$ annihilating $F_2$. (The motivation for the notation $\fgl(U)_{0,*}$ is explained in \S\ref{oscfiltsect}.)  Since $\fr$ is reductive, we may decompose $\fgl(U)_{0,\ast}=\fr\op\fr\upperp$ as an $\fr$-module.
In the notation of \S\ref{edssect},
\begin{displaymath}
  V \ \simeq \ L^*\ot T \, , \quad 
  W \ \simeq \ (L^*\ot N)\op (T^*\ot N) \, , \quad 
  A \ \simeq \ \fr\upperp \, , 
\end{displaymath}
with $L^*\ot N$ in the first derived system.
That  $\fr\upperp \subset V^*\ot W$    may be seen as follows
\begin{eqnarray}
  \nonumber 
  d(\o_{T^*\ot N}-F_2(\o_{L^*\ot T})) & = & 
  - \ \o_{T^*\ot L}\ww\o_{L^*\ot N} \ - \ \o_{T^*\ot T}\ww\o_{T^*\ot N} \\
  & & \nonumber
  - \ \o_{T^*\ot N}\ww\o_{N^*\ot N} \ 
  + \ F_2 ( \o_{L^*\ot L}\ww\o_{L^*\ot T} ) \\
  & & \nonumber
  + \ F_2 ( \o_{L^*\ot T}\ww\o_{T^*\ot T} \ 
  + \ \o_{L^*\ot N}\ww\o_{N^*\ot T} ) \\
  & \equiv & \label{f2calc}
  - \o_{T^*\ot T}\ww F_2(\o_{L^*\ot T}) \ 
  - \ F_2(\o_{L^*\ot T})\ww\o_{N^*\ot N} \\ \nonumber
  & & - F_2(-\o_{L^*\ot L}\ww\o_{L^*\ot T} \ 
      - \ \o_{L^*\ot T}\ww\o_{T^*\ot T}) \ \tmod I  \\ \nonumber
 & \equiv &  (\o_{0,\ast} \, . \, F_2)\ww\o_{L^*\ot T} \ \tmod I \\ \nonumber
 & \equiv & (\o_{\fr\upperp} \, . \, F_2)\ww\o_{L^*\ot T} \ \tmod I \, .
\end{eqnarray}
Above, $\w_{0,\ast} \, . \, F_2$ denotes the action of the $\fgl(U)_{0,\ast}$--valued component $\w_{0,\ast}$ of the Maurer-Cartan form on $F_2 \in S^2T^* \ot N$.  Recall that $\fr$ is the annihilator of this action.    By definition $\w_{0,\ast} \, . \, F_2 = (\w_\fr + \w_{\fr^\perp}) \, . \, F_2= \w_{\fr^\perp} \, . \, F_2$.

\subsection{Third order Fubini systems}\label{sec:Fub3}
Fix    $F_2\in L \ot S^2T^*\ot N$ and
$F_3\in L^2 \ot S^3T^*\ot N$.  
The third order Fubini system corresponding to $(F_2,F_3)$ is
\be\label{f3sys}
  \begin{array}{r@{ \ = \ }l}
I_{\tFub_3} & \{\o_{L^*\ot N} \, , \ \o_{T^*\ot N}-F_2(\o_{L^*\ot T}) \, , \
 \o_{\fr\upperp}.F_2-F_3(\o_{L^*\ot T})\} \, \\
 J_{\tFub_3} & \{I_{\tFub_3} \, , \ \o_{L^*\ot T}\}.
\end{array}\ene
where the last term in $I_{\tFub_3}$ is $N\ot T^*\ot T^*\ot L$-valued,
and   $\o_{\fr\upperp}.F_2$ is as in \S\ref{sec:Fub2}.  Here
\begin{displaymath}
  V \ \simeq \  L^*\ot T \, , \quad
  W \ \simeq \ (L^*\ot N)\op (T^*\ot N)\op \fr\upperp \, , 
\end{displaymath}
with the first two terms of $W$ in the first derived system.  The 
tableau $A\simeq (L\ot T^*)\op (T\ot N^*)$ sits in $W\ot V^*$ by a 
calculation similar to \eqref{f2calc}.

Integral manifolds of the third order Fubini system are adapted frame
bundles $\cF^3 \subset GL(U)$ over submanifolds of projective space whose second fundamental form is isomorphic to $F_2$ at each point and whose Fubini cubic is normalizable
to $F_3$ at each point.  
 Note that these systems admit reduced prolongations.  The reductions (or normalizations) are by the subalgebra
$\fn\subset ( L \ot  T^* ) \op ( T \ot N^* ) \op ( L \ot N^* )$
such that $\fn.(\fr\upperp . \, F_2 - F_3(L^*\ot T))=0$.

In this paper we will be concerned with $(F_2,F_3)$ modeled on
the Fubini forms of a homogeneous variety.

\section{A bi-grading of $\fgl(U)$}
\subsection{The osculating filtration}\label{oscfiltsect}

Given a submanifold $X\subset \BP U$, and $x\in X$, the {\it osculating filtration at $x$}
$$U_0\subset U_1\subset\cdots \subset U_\ell=U$$
is defined by 
$$
U_0=\hat x,\ U_1=\hat T_xX,\ U_2=U_1+F_2(L^* \ot S^2T_xX)\hd \ U_r=U_{\ell-1}+F_\ell(L^*{}^{\ell-1}\ot S^\ell T_xX).
$$   The $F_j$ are the Fubini forms of \S\ref{fubformssect}.  We may reduce the frame bundle $\cF^1_X$ to framings adapted to the osculating sequence by specifying 
$e = (e_0 , e_\a , e_{\m_2} , \ldots , e_{\m_\ell}) \in \cF^1_X$ so that  
$U_k = \tspan \{ e_0 , e_\a , e_{\m_2} , \ldots , e_{\m_k} \}$.  (The indices $\a$ and $\m_k$ 
respectively range over $1 , \ldots , n$ and $\tdim U_{k-1} + 1 , \ldots , \tdim U_k$.)  From now on we work on this reduced frame-bundle, denoted $\cF^\ell_X \subset \cF^1_X$.   

At each point of $\cF^\ell_X$ we obtain a splitting of $U$.  This induces a splitting
$$
\fgl(U)=\oplus\fgl(U)_{k,\ast} \, .
$$
The asterisk above is a place holder for a second  splitting given (in \S\ref{rootgradingsect}) by the representation theory when $X = G/P$. 

The osculating filtration of $U$   determines a refinement of the Fubini forms.  Let $N_k=U_k/U_{k-1}$ and decompose $F_s=\op_kF_{k,s}$ so that $F_{k,s}: L^*{}^{s-1}\ot S^sT_xX \to N_k$. 
Although the Fubini forms do not descend to well-defined tensors on $X$, the {\it fundamental forms} $F_{k,k}$ do.  By definition, $F_{k,k} : L^*{}^{k-1} \ot S^kT_xX \to N_{k,x}X$ is surjective.     
\medskip

\noindent {\it Remark.}
The only nonzero Fubini forms of a homogeneously embedded CHSS are the fundamental forms.  For the adjoint varieties, the only nonzero
Fubini forms are $F_{2,2},F_{2,3},F_{2,4}$.
\medskip

\subsection{\boldmath The root grading and the $(I_p,J_p)$ systems \unboldmath}\label{rootgradingsect}

Let $\tilde\fg$ be a complex semi-simple Lie algebra with a fixed set of simple roots $\{ \a_1\hd \a_r \}$, and corresponding fundamental weights $\{\w_1 , \ldots , \w_r \}$.  Let $I \subset \{1 , \ldots , r\}$, and consider the irreducible representation $\mu:\tilde\fg\ra\fgl(U)$ of highest weight $\l=\sum_{i \in I} \l^{i}\o_{i}$.
Set $\fg=\mu(\tilde\fg)$, and let $\mu (G)\subset GL(U)$
be the associated Lie group so that
$ G/P \subset \BP U$ is the orbit of a highest
weight line.  Write $P=P_I\subset G$ for the parabolic subgroup obtained by deleting negative root spaces corresponding to roots having a nonzero coefficient on any of the simple roots $\a_{i}$, $i \in I$.

Since $\tilde\fg$ is reductive, we have a splitting $\fsl(U)=\fg\op \fg\upperp$, where $\fg\upperp$ is the $\tilde\fg$-submodule of $\fgl(U)$ complementary to $\fg$. 
Let $\o\in \O^1(GL(U),\fgl(U))$ denote the Maurer-Cartan form of $GL(U)$, and let $\o_{\fg}$ and $\o_{\fg^\perp}$ denote the components of $\o$ taking values in $\fg$ and $\fg^\perp$, respectively.  

The bundle  $\cF^\ell_{G/P}$ 
admits a reduction to a bundle $\cF^G_{G/P} =\mu (G)$.  On this bundle  the Maurer-Cartan form pulls-back to take values in $\fg$,
that is,
$\w_{\fg^\perp} = 0$.  Conversely, all $\tdim (G)$ 
dimensional integral
manifolds of the system $I=\{\w_{\fg^\perp}\}$
are conjugates of $\mu (G)$.  That is, the $\w_{\fg^\perp} = 0$ system is rigid.  Theorem \ref{transthm} establishes the rigidity of   weaker systems, subject to the partial vanishing of
components of a Lie algebra cohomology group.
 
Let $Z=Z_I\subset \ft$ be the {\it grading element}
corresponding to $\sum_{i\in I}\a_{i}$. The grading element $Z_i$ for a simple root $\a_{i}$ has the  property that $Z_i(\a_j)=\d^i_j$. In general $Z=\sum_{i \in I} Z_{i}$. Thus, if   $(c^{-1})$ 
denotes the   inverse of the Cartan matrix, then given a weight $\nu=\sum\nu^j\o_j$,  
\be\label{charelteqn}
Z(\nu)=\sum_{{1\le j \le r}\atop{i \in I}}\nu^j(c\inv)_{j,i} \, .
\ene
The grading element induces a $\bZ$-grading of $\fg = \oplus_{-k}^k \fg_k$.  To determine $k$ in the case $P=P_{\a_i}$ is a maximal parabolic, let $\tilde\a$ denote the highest
root. Given $\tilde\a=\sum m_j\a_j$, we have $k=m_i$.

The module $U$ inherits a $\BZ$-grading 
$$
U=U_{Z(\l)}\op U_{Z(\l)-1}\op\cdots\op U_{Z(\l)-f} \, .
$$
The $U_j$ are eigen-spaces for $Z$.
This grading is compatible with the action of $\tilde \fg$: $\mu(\tilde\fg_i).U_j\subset U_{i+j}$. 
We   adopt the notational convention   of shifting the grading
on $U$  to begin
at zero. 
The component $U_0$ (formally named $U_{Z(\l)}$) is one dimensional,
corresponding to the highest weight line of $U$, and $G \cdot \bP U_0 = G/P \subset \bP U$. (The 
labeling of the grading on $\fgl(U)=U^*\ot U$
is independent of our shift convention.)

Note, in particular, that the vector space 
$\hat T_{[\tId]}(G/P)/\hat{\tId}\simeq \fg/\fp$
is graded from $-1$ to $-k$.  The osculating grading on $U$ induces gradings of
$\fgl(U)$, $\fg$ and $\fg^\perp$.


 
We write
$$
\fgl(U)=\bigoplus_{s,j} \fgl(U)_{s,j}
$$
where the first index refers to the osculating grading (\S\ref{oscfiltsect})
induced by $G/P\subset \BP U$ and the second the root grading. We adopt the notational convention
$$
\fgl(U)_j=\bigoplus_{s } \fgl(U)_{s,j} \, ;
$$
so
if there is only one index, it refers to the root grading. Note that although $Z(\l)$ need not be an integer, $f$ is.  So the grading of $\fgl(U)$ is indexed by integers $-f\hd f$.  Moreover, this grading of $\fgl(U)$ is independent of shift convention.
\medskip

We began with frame bundles $\cF \subset \tGL(U)$.  But rescaling $e_0$ we may assume that $\cF \subset \tSL(U)$, and we will do so from now on.  We define the $(I_{p},J_{p})$ system on
$SL(U)$ by
$I_p=\{ \o_{\gp{\leq p}}\}$, $J_p=\{ I_p, \o_{\fg_-}\}$.
The EDS we will study are the filtered systems $(I_p^\textsf{f},\Omega)$,
which are weaker than the $(I_p,J_p)$ systems.  (See \S\ref{filtersect}.)

\medskip

\noindent {\it Remark.} The osculating grading coincides
(up to a change of sign)
with the Lie algebra grading if and only if $G/P$ is CHSS.
\medskip

  
\section{The adjoint varieties}\label{motexsecta}
In this section we  describe the 
fundamental adjoint representations from a uniform
perspective and prove that integral
manifolds of  the third order Fubini
system for them are automatically integral
manifolds for the corresponding $(I_{-1},J_{-1})$-system.

\smallskip

\noindent {\it Remark.}  We have the following geometric models for   adjoint varieties $Z^G_\tad$.
\begin{list}{$\circ$}
  {\usecounter{cl}
   \setlength{\leftmargin}{15pt}
   \setlength{\labelwidth}{10pt}
   }  
\item 
${A_n} = SL_{n+1}\bC$:  The flag variety $\BF_{1,n}(\BC^{n+1})$ 
of lines in hyperplanes in $\BC^{n+1}$.  
\item 
${C_n} = Sp_{2n}\bC$:  The Veronese variety $v_2(\pp{2n-1})\subset \BP (S^2\BC^{2n})$.
\item 
$G = SO(n)$: The variety $G_Q(2,n)$ of $2$-planes in $\BC^n$ on
which the quadratic form restricts to be zero.
\item 
$G = {G_2}$: The variety $G_{\mathrm{null}}(2,\mathrm{Im}(\BO))$ of $2$-planes in the imaginary complexified octonions to which the octonionic multiplication restricts
to be identically zero.
\item $G = {E_6}$: The variety of $\pp 6$'s on
the Cayley plane.
\end{list}
\smallskip

  The highest root $\tilde\a$ induces a five-step grading $\fg=\fg_{-2}\op\fg_{-1}\op \fg_0\op \fg_1\op\fg_2$.  (The corresponding grading of $U=\tilde\fg$ is $U_j = \fg_{j-2}$.)  The adjoint variety is the $G$--orbit of $\bP\fg_2$. 

Fix a Chevalley basis of $\fg$.  Decompose $\fg_0 = \ff + \BC\{Z_{\tilde\a}\}$ into the semi-simple Levi factor $\ff$, and the one-dimensional (assuming $\fg\neq \fa_n$) $\BC\{Z_{\tilde\a}\}$ .  In order to obtain the representation $\fg=\tad(\tilde \fg) \subset \fgl(U)$ we make the following observations.  First, $ Z_{\tilde\a} $ acts on $\fg_j$ by $j \, \tId$; and $\ff$ acts trivially on $\fg_{\pm2}$.
Let $\tii: \ff\ra \fgl(\fg_{-1}) = \fg_{-1}\ot\fg_{-1}^*$ denote the action of 
$\ff$ on $\fg_{-1}$.  (Our choice of notation is made to be compatible with the geometry that enters later.)  The Killing form $B$ on $\fg$ allows us to identify $\fg_{\pm1} \simeq \fg_{\mp1}{}^*$, so that we may regard $\tii$ as an element of $\ff^*\ot \fg_{-1}\ot \fg_1$.  The bracket $[\cdot,\cdot]: \fg_{-1}\times\fg_{-1} \to \fg_{-2}$ induces a symplectic form $\Omega \in \La 2\fg_{-1}^*$.  The linear map $\fg_{-1}\ra\fg_{-1}^*=\fg_1$ induced by $\Omega$ is an $\ff$--module isomorphism. (That is, the representation $\tii$ is symplectic.)
Let $\mathit{II}^{-2}=\mathit{contr}(B\ot \tii\ot \O)\in \fg_{-1}^*\ot \fg_{-1}^*\ot \ff$
denote the natural contraction
$$
(\ff \ot \ff )\ot (\ff^*\ot  \fg_{-1}^*\ot \fg_{-1})
\ot (\fg_{-1}^*\ot \fg_{-1}^*)\ra
\ff\ot \fg_{-1}^*\ot \fg_{-1}^* \, .
$$
Because the representation $\tii$ is symplectic,
the image is in $\ff\ot S^2\fg_{-1}^*$.

The  canonical identification $\fg_{\pm1}\simeq \fg_{\mp1}^*$ allows us to   treat $\mathit{II}^{-2}$, $\tii$ and $\O$ as elements of $\ff\ot S^2\fg_1$, $\ff^*\ot \fg_1^*\ot \fg_1$ and $\La 2\fg_1$, respectively.

The adjoint variety admits a frame bundle $\mu(G) = \cF^G_\tad\subset\tGL(U)$ on which the Maurer-Cartan form pulls-back to take values in $ \fg\subset \fgl(U)$.  Decomposed with respect to the bi-grading, $\w$ takes the form
\be\begin{pmatrix}
2  Z_{\tilde\a}  & \fg_{1} &\fg_2 & 0 &0 &0\\
\fg_{-1} &  Z_{\tilde\a} \, \tId +\tii (\ff) &
\O(\fg_1) & \tii(\O(\fg_1)) & -\tfrac12 \, \fg_2 \, \tId  &0\\
\fg_{-2} &-\frac 12\O(\fg_{-1}) &0 &0 &-\frac12 \, \fg_1
& -\tfrac12 \, \fg_2 \\
0&\mathit{II}^{-2}(\fg_{-1})& 0 & \tad(\ff) & -\mathit{II}^{-2}(\O(\fg_1))
&0\\
0&-\fg_{-2} \, \tId & -\fg_{-1}&
\tii(\fg_{-1})&- Z_{\tilde\a} \tId + \tii(\ff)
&\O(\fg_1)\\
0&0&-2 \, \fg_{-2}&0&-\O(\fg_{-1}) &-2 Z_{\tilde\a}
\end{pmatrix} \, .
\ene

The last non-zero fundamental form is the second, so the osculating filtration has length two.  The tangent space $T = N_1$ and the normal space $N = N_2$ decompose as $N_{1,-1} \op N_{1,-2}$ and $N_{2,-2} \op N_{2,-3} \op N_{2,-4}$, respectively, under the root grading.  We write $T_j = N_{1,j}$ and $N_{j} = N_{2,j}$.
Notice that $\mathit{II}^{-2}$ is the restriction to $N_{-2}^*$ of $\mathit{II}: N^* \to S^2T^*$.
  The induced bi-grading on $\fgl(U)_{\textrm{osc},\textrm{alg}}$
is indicated in the table below.
\begin{center}\renewcommand{\arraystretch}{1.3}
\begin{tabular}{c||c|c|c|c|c|c}
    & $L^*$ & $T_{-1}{}^*$ & $T_{-2}{}^*$ & $N_{-2}{}^*$ & $N_{-3}{}^*$ 
    & $N_{-4}{}^*$ \\  \hline \hline
 $L$ & 
 (0,0) & (-1,1) & (-1,2) & (-2,2) & (-2,3) & (-2,4) \\ \hline
 $T_{-1}$ & 
 (1,-1) & (0,0) & (0,1) & (-1,1) & (-1,2) & (-1,3) \\ \hline
 $T_{-2}$ & 
 (1,-2) & (0,-1) & (0,0) & (-1,0) & (-1,1) & (-1,2) \\ \hline
 $N_{-2}$ & 
 (2,-2) & (1,-1) & (1,0) & (0,0) & (0,1) & (0,2) \\ \hline
 $N_{-3}$ & 
 (2,-3) & (1,-2) & (1,-1) & (0,-1) & (0,0) & (0,1) \\ \hline
 $N_{-4}$ & 
 (2,-4) & (1,-3) & (1,-2) & (0,-2) & (0,-1) & (0,0)
\end{tabular}
\end{center}

\begin{proposition}\label{prop:Fub3=>-1}
Every integral manifold of the third-order Fubini system $(I_{\tFub_3},J_{\tFub_3})$ for
a given adjoint variety is an integral
manifold of the $(I_{-1},J_{-1})$ system for the
same adjoint variety.
\end{proposition}

\begin{proof}
Suppose that $\cF \subset SL(U)$ is an integral manifold of the third-order Fubini system.    We wish to show that the $\fg^\perp_{*,<0}$--valued component of the Maurer-Cartan form vanishes when pulled-back to $\cF$.  That the $\fg^\perp_{>0,*}$--valued component vanishes is an immediate consequence of the injectivity of the second fundamental form $F_2$ on each homogeneous component.

Referring to the table above, we see that there remain
four blocks of the component of the Maurer-Cartan form 
in $\gp{*,<0}$ to consider: the three $(0,-1)$ blocks $\w_{T_{-2}\ot T_{-1}^*}$, $\w_{N_{-3}\ot N_{-2}^*}$ and $\w_{N_{-4}\ot N_{-3}^*}$; and the singleton $(0,-2)$ block $\w_{N_{-4}\ot N_{-2}^*}$.  The third Fubini form is defined by (3.5) of \cite[\S3.5]{EDSBK}.  The vanishing of the $\fg^\perp$--component of the first two blocks is a consequence of the $S^3 T_{-1}^* \ot N_{-3}$ component of $F_3$.  (This is the only nonzero component of $F_3$.)  The vanishing of the $\fg^\perp$--component of the third and fourth blocks is given by the $S^3(T^*) \ot N_{-4}$--component of $F_3$.
\end{proof}

\section{Filtered EDS}\label{filtersect}

In \cite{HY}, it was observed (in different language) that
for rigidity problems associated to CHSS, the  $(I_{-1},J_{-1})$ system could be
proved to be rigid using Lie algebra cohomology. More precisely,  the
Spencer differential coincides with the Lie algebra cohomology differential $\partial^1_1$
and the admissible normalizations coincide with the image of the Lie algebra cohomology map $\partial^0_1$. This correspondence breaks down when $k>1$ (see Remark \ref{breakdownrem} below), but it can be restored
with the use of filtered EDS and simultaneously studying the prolongation and torsion.

\begin{definition} Let $\Sigma$ be a manifold equipped with a filtration of its
tangent bundle
$T^{-1}\subset T^{-2}\subset \cdots \subset T^{-f}=T\Sigma$.
Define an \emph{$r$-filtered} Pfaffian EDS on $\Sigma$ to be a filtered ideal
$I\subset T^*\Sigma$ whose
integral manifolds are the immersed submanifolds $i: M\ra \Sigma$
such that 
$i^*(I_{u})|_{i^*(T^{u-r})}=0$ for all $u$, with the convention that
$T^{-s}=T\Sigma$ when $-s<-f$.
\end{definition}

Another way to view filtered EDS is to consider an ordinary EDS on the
total space of the sum of the bundles
$I_u\otimes (T\Sigma/T^{u+r})$. In our case these bundles will be trivial with
fixed vector spaces as models.

Define $(I^\textsf{f}_{p},\Omega)$ to be the $(p+1)$-filtered EDS on $GL(U)$ with filtered ideal
  $I^\textsf{f}_{p}:=\o_{\fg^{\perp}_{\leq p}}$ and independence condition
$\Omega$ given by the wedge product of the forms in $\o_{\fg_{-}}$.   We may view this as an ordinary
EDS on 
$$GL(U)\times \left( [\fg^{\perp}_{p}\ot (\fg_{-2}\op \cdots \op \fg_{-k})^*]
\oplus [\fg^{\perp}_{p-1}\ot (\fg_{-3}\op \cdots \op \fg_{-k})^*]
\oplus \cdots \oplus [\fg^{\perp}_{p-k+2}\ot  \fg_{-k}^*] \right)
$$
where, giving $\fg^{\perp}_i\ot \fg_{-j}^*$ linear coordinates $\l_{i,j}$, we
have
\begin{eqnarray} \label{psystem} 
  I^\textsf{f}_p & = &
  \big\{ \ \o_{\gp s} \, , \ s\leq p-k+1 \, ; \quad
  \o_{\gp{p-k+2}}-\l_{p-k+2,k}(\o_{\fg_{-k}}) \, , \\
  & & \nonumber \hspace{10pt}
  \o_{\gp{p-k+3}}-\l_{p-k+3,k}(\o_{\fg_{-k}})-\l_{p-k+3,k-1}(\o_{\fg_{-k+1}}) \, , \ \ldots \\
  & & \nonumber \hspace{10pt}
  \o_{\gp{p}}-\l_{p,k}(\o_{\fg_{-k}})-\cdots -\l_{p,2}(\o_{\fg_{-2}}) \ \big\}
\end{eqnarray}
\smallskip

\noindent
We henceforth assume {\bf \boldmath $p\geq -1$\unboldmath}.

\subsection{The $p=-1$, $k=2$ case}
\label{sec:warmup}
We give a proof of the main result in \S\ref{sec:gen_case} below, but for the reader's convienence we work out the $p=-1$ and $k=2$ case explicitly to illustrate the central ideas.

Abbreviate
$$ \w_{\fsl(U)_s} \ =: \ \w_s.$$
Notice that $\fg_{s}^\perp = \fsl(U)_{s}$ for all $s \le -3$.  
so that $\w_{\fg^\perp_s} = \w_{s}$, for all $s \le -3$.  Thus the ideal $I^\textsf{f}_{-1}$ is 
\begin{displaymath}
  I^\textsf{f}_{-1} \ = \ \left\{ 
    \w_{\fg^\perp_{-1}} - \lambda_{-1,2}(\w_{\fg_{-2}}) \, , \ 
    \w_{\fg^\perp_{-2}} \, , \ 
    \w_{-3} \, , \ \ldots \, , \ \w_{-f} 
  \right\} \, .
\end{displaymath}

The calculations that follow utilize the Maurer-Cartan equation (\ref{eqn:mc}), and
the facts that $[\fg,\fg] \subset \fg$ and $[\fg,\fg^\perp] \subset \fg^\perp$.  It is easy to see that $\td \w_{s} \equiv 0$ modulo $I^\textsf{f}_{-1}$ when $s \le -4$.  Next, computing modulo $I^\textsf{f}_{-1}$, 
\begin{eqnarray}
  - \td \, \w_{-3} & \equiv & \label{eqn:-3}
  \left[ \w_{\fg_{-2}} , \w_{\fg^\perp_{-1}} \right] 
  \ \equiv \
  \left[ \w_{\fg_{-2}} , \lambda_{-1,2}(\w_{\fg_{-2}}) \right]  \, , \\
  -\td \, \o_{\gp{-2}} & \equiv & \label{eqn:-2}
  \left[\o_{\fg_{-2}} , \o_{\gp 0} \right] \ + \ 
  \left[\o_{\fg_{-1}},\w_{\gp {-1}} \right] \ + \
  \left[ \w_{\gp {-1}}, \w_{\gp {-1}} \right]_{\fg\upperp}\\
  & \equiv & \nonumber
  \left[\o_{\fg_{-2}} , \o_{\gp 0} \right] \ + \ 
  \left[\o_{\fg_{-1}},\lambda_{-1,2}(\w_{\fg_{-2}}) \right] \ + \
  \left[\lambda_{-1,2}(\w_{\fg_{-2}}) , 
        \lambda_{-1,2}(\w_{\fg_{-2}}) \right]_{\fg\upperp} \, .
\end{eqnarray}
Finally, 
\begin{equation}\label{eqn:-1}
\renewcommand{\arraystretch}{1.3}
\begin{array}{l}
  -\td \left(\o_{\gp{-1}}-\l_{-1,2}(  \o_{\fg_{-2}})\right)  \ \equiv \\
  \hspace{80pt} 
  \left[ \w_{\fg_{-2}} , \w_{\fg^\perp_1} \right] \ + \ 
  \left[ \w_{\fg_{-1}} , \w_{\fg^\perp_0} \right] \ + \ 
  \left[ \w_{\fg_{-1}^\perp} , \w_{\fg_0} \right] \ + \ 
  \left[ \w_{\fg^\perp_{-1}} , \w_{\fg^\perp_{0}} \right]_{\fg^\perp} \\
  \hspace{80pt} + \ 
  \td \lambda_{-1,2} ( \wedge \w_{\fg_{-2}} ) \ - \
  \lambda_{-1,2} \left(
    \left[ \w_{\fg_{-2}} , \w_{\fg_0} \right] + 
    \left[ \w_{\fg_{_1}} , \w_{\fg_{-1}} \right] +
    \left[ \w_{\fg^\perp_{-1}} , \w_{\fg^\perp_{-1}} \right]_{\fg}
  \right) \\
  \hspace{65pt} \equiv  \
  \left[ \w_{\fg_{-2}} , \w_{\fg^\perp_1} \right] \ + \ 
  \left[ \w_{\fg_{-1}} , \w_{\fg^\perp_0} \right] \ + \ 
  \left[ \lambda_{-1,2}(\w_{\fg_{-2}}) , \w_{\fg_0} \right] \\
  \hspace{80pt} + \
  \left[ \lambda_{-1,2}(\w_{\fg_{-2}}) , 
         \w_{\fg^\perp_{0}} \right]_{\fg^\perp} \ + \ 
  \td \lambda_{-1,2} ( \wedge \w_{\fg_{-2}} )  \\
  \hspace{80pt} - \
  \lambda_{-1,2} \left(
    \left[ \w_{\fg_{-2}} , \w_{\fg_0} \right] + 
    \left[ \w_{\fg_{_1}} , \w_{\fg_{-1}} \right] +
    \left[ \lambda_{-1,2}(\w_{\fg_{-2}}) , \lambda_{-1,2}(\w_{\fg_{-2}}) \right]_{\fg}
  \right)
\end{array}
\end{equation}
Here the bracket is as indicated in \S\ref{sec:mc} and $[\cdot , \cdot]_{\fg}$ (resp. $[\cdot,\cdot]_{\fg\upperp}$) denotes
the component of the bracket taking values in $\fg$ (resp. $\fg\upperp$).

The three differentials above must vanish on an integral element  modulo $I^\textsf{f}_{-1}$.  Notice that the vanishing of (\ref{eqn:-2}) implies 
$$
\o_{\gp 0}= \l_{0,1}(\o_{\fg_{-1}})+ \l_{0,2}(\o_{\fg_{-2}})
$$
for some $\l_{0,j}\in \gp 0\ot \fg_{-j}^*$.

Now consider the degree one ($1=p+2$) homogeneous component of (\ref{eqn:-3},\ref{eqn:-2},\ref{eqn:-1}),
simultaneously examining the torsion and tableau.
In order the have an integral element, ${\mathbf\l}_1:=\oplus_{s=-1}^0\l_{s,1-s}$ must be in the kernel of the map
$$
\d_1:
\oplus_{s=-1}^0(\gp{s}\ot \fg^*_{-s-1})
\ra   (\gp{-1}\ot  \fg^*_{-1}\ww \fg^*_{-1})
\op (\gp{-2}\ot \fg^*_{-1}\ww \fg^*_{-2})
$$
defined as follows.  Given $u_{-1},v_{-1}\in \fg_{-1}$, 
\begin{equation}\label{d1der}
\d_1({\mathbf\l}_1)(u_{-1}\ww v_{-1})
=[\l_{0,1}(u_{-1}),v_{-1}]+[u_{-1},\l_{0,1}(v_{-1})]-\l_{-1,2}([u_{-1},v_{-1}]).
\end{equation}
For $u_{-1}\in \fg_{-1},v_{-2}\in \fg_{-2}$
$$
\d_1({\mathbf\l}_1)(u_{-1}\ww v_{-2})
=[\l_{0,1}(u_{-1}),v_{-2}]+[u_{-1},\l_{-1,2}(v_{-2})] \, .
$$
That is,
$\d_1=\partial^1_1$, where $\partial^1_1$ is the Lie algebra cohomology differential 
described in \S\ref{cohrulessect}.

\begin{remark}\label{breakdownrem}
Note that had we instead used the $(I_{-1},J_{-1})$ system in equation \eqref{d1der}, the $\l_{-1,2}$ term
would be missing and we would not have   $\partial^1_1$.
The map $\delta_1$ is akin to an `augmented Spencer differential' it addresses both the torsion and prolongation in homogeneous degree 1.
\end{remark}

Moreover, $\fg^\perp_{1}  = \fn\cap  \fgl(U)_{ 1}$ (cf. \S\ref{redprolongsect}), and the Lie algebra cohomology denominator $\partial^0_{1}(\fg^\perp_{1})$ is the space of admissible normalizations of the prolongation coefficients ${\mathbf \l}_{1}$.  Thus, the vanishing of $H^1_{1}(\fg_-,\fg^\perp)$ implies that normalized integral manifolds of the $(I^f_{-1},\Omega)$ system are in one to one
correspondence with integral manifolds of the $(I^f_{0},\Omega)$ system.

\subsection{The general case}
\label{sec:gen_case}

Define 
$$ \lambda_{s,*} \ := \ \bigoplus_{j=p+2-s}^k \lambda_{s,j} \ \in \ 
   \bigoplus_{j=p+2-s}^k  \fg^\perp_s \ot \fg_{-j}^* \, ,
$$
and set 
$$ \lambda_{s,*} \left( \w_{\fg_{<0}} \right) \ := \ 
   \sum_{j=p+2-s}^k \lambda_{s,j} \left( \w_{\fg_{-j}} \right) \, .
$$
Defining $\lambda_{s,j}=0$ when $s+j<p+2$ allows us to write the generators of the ideal $I^\textsf{f}_p$ compactly as 
$\w_{\fg^\perp_s} - \lambda_{s,*}(\w_{\fg_{<0}})$.  

Observe that, if $s < -k$, then $\fg^\perp_s = \fsl(U)_s$, so that $\w_s \equiv 0$ modulo $I^\textsf{f}_p$.  Also, $\w_{\fg^\perp_{-k}} \equiv 0$ modulo $I^\textsf{f}_p$.  Now, the Maurer-Cartan equation (\ref{eqn:mc}) yields, modulo $I^\textsf{f}_p$,  
\begin{eqnarray}
  & & \nonumber
  -2 \, \td (\w_{\fg^\perp_s} - \lambda_{s,*}(\w_{\fg_{<0}}))  \\ 
  & & \nonumber
  \qquad \equiv \ 
  \sum_{t+t'=s} \left[ \w_t , \w_{t'} \right]_{\fg^\perp} \ 
  + \ 2 \, \td \lambda_{s,*} \wedge \w_{\fg_{<0}} \ 
  + \ \lambda_{s,*} 
      \left( \sum_{-k\le t+t'<-1} \left[ \w_t , \w_{t'} \right]_{\fg} \right)  \\
  & & \qquad \equiv \   \nonumber
  \sum_{t=-k}^{s-p-1} \left[ \w_{\fg_t} + \lambda_{t,*}(\w_{\fg_{<0}}) , 
                             \w_{\fg_{s-t}}+\w_{\fg^\perp_{s-t}} \right]_{\fg^\perp} \\
  & & \qquad\qquad + \   \label{eqn:egads}
  \sum_{t=s-p}^p \left[ \w_{\fg_t} + \lambda_{t,*}(\w_{\fg_{<0}}) , 
                        \w_{\fg_{s-t}} + \lambda_{s-t,*}(\w_{\fg_{<0}})  \right]_{\fg^\perp} \\
  & & \qquad\qquad + \   \nonumber
  \sum_{t=p+1}^{s+k} \left[ \w_{\fg_t}+\w_{\fg^\perp_t} , 
                            \w_{\fg_{s-t}} + \lambda_{s-t,*}(\w_{\fg_{<0}})  \right]_{\fg^\perp}
   \qquad + \ 2 \, \td \lambda_{s,*} \wedge \w_{\fg_{<0}} \\
  & & \qquad \qquad      \nonumber
  - \sum_{j=p+2-s}^k \lambda_{s,j} 
      \Bigg(
      \sum_{t=-k}^{-j-p-1} \left[ \w_{\fg_t} + \lambda_{t,*}(\w_{\fg_{<0}}) , 
                             \w_{\fg_{-j-t}}+\w_{\fg^\perp_{-j-t}} \right]_{\fg} \\
  & & \hspace{130pt} \ + \   \nonumber
      \sum_{t=-j-p}^p \left[ \w_{\fg_t} + \lambda_{t,*}(\w_{\fg_{<0}}) , 
                        \w_{\fg_{-j-t}} + \lambda_{-j-t,*}(\w_{\fg_{<0}}) \right]_{\fg} \\
  & & \hspace{130pt} \ + \   \nonumber
      \sum_{t=p+1}^{k-j} \left[ \w_{\fg_t}+\w_{\fg^\perp_t} , 
                            \w_{\fg_{-j-t}} + \lambda_{-j-t,*}(\w_{\fg_{<0}}) \right]_{\fg}
      \Bigg) \, .
\end{eqnarray}
Observe that, for $s \le -2k$, $\td \w_{\fg_s^\perp} \equiv 0$ modulo $I^\textsf{f}_p$; that is, for $s\le-2k$, $\w_{\fg_s^\perp}$ lies in the first derived system $I'$.  


The necessary vanishing of (\ref{eqn:egads}) on integral elements tells us that there exist $\lambda_{s,j}$, taking values in $\fg^\perp_s \ot \fg_{-j}{}^*$, such that 
\begin{equation}\label{eqn:prolong}
  \w_{\fg^\perp_s} = \sum_{j=1}^k \lambda_{s,j} (\w_{\fg_{-j}}) \ , \quad 
  p+1 \le s \le p+k \, .
\end{equation}  
The vanishing of (\ref{eqn:egads}) also imposes consraints on the $\w_{\fg_{s}}$ and $\td \lambda_{s,*}$, but these terms will not appear when we consider the first non-zero component of the augmented Spencer differential $\delta_{p+2}$ below.  (Of course,  the $\w_{\fg_s}$, $s\ge0$, are dual to the Cauchy characteristics, and thus ought not influence the computation.)   

Careful consideration of (\ref{eqn:egads}) shows that the right-hand side contains no components of homogeneous degree less than $p+2$.  We identify the necessary condition that the degree $p+2$ homogeneous component of the right-hand side of (\ref{eqn:egads}) vanish on an integral element by evaluating the expression on $u \in T^{-i}/T^{1-i} \simeq \fg_{-i}$, $v \in T^{-j}/T^{1-j} \simeq \fg_{-j}$ with $s=(p+2)-(i+j)$, $1 \le i,j \le k$.  Utilizing (\ref{eqn:prolong}), we have 
\begin{eqnarray}
  & & \nonumber
  -2 \, \td (\w_{\fg^\perp_s} - \lambda_{s,*}(\w_{\fg_{<0}}))(u,v)  \\ 
  & & \qquad \equiv \   \nonumber
  \hspace{20pt}
  \sum_{t=-k}^{s-p-1} \left[ \w_{\fg_t}(u) , 
                             \lambda_{s-t,*}( \w_{\fg_{<0}}(v)) \right] \\
  & & \qquad\qquad + \   \nonumber
  \sum_{t=s-p}^p 
         \left[ \w_{\fg_t}(u) , \lambda_{s-t,*} ( \w_{\fg_{<0}}(v) ) \right] \ 
     + \ \left[ \lambda_{t,*} ( \w_{\fg_{<0}}(u)) , \w_{\fg_{s-t}}(v)  \right] \\
  & & \qquad\qquad + \   \nonumber
  \sum_{t=p+1}^{s+k} \left[ \lambda_{t,*} ( \w_{\fg_{<0}}(u)) , 
                            \w_{\fg_{s-t}}(v) \right]  \\
  & & \qquad \qquad      \nonumber
  - \ \lambda_{s,*} \left( \left[ \w(u) , \w(v) \right] \right) \\ 
  & & \qquad  = \   \nonumber
  \hspace{20pt}
  \left[ \w_\fg (u) , \lambda_{s-i,*} ( \w_{\fg_{<0}}(v)) \right] 
  \ - \ 
  \left[ \w_\fg (v) , \lambda_{s-j,*} ( \w_{\fg_{<0}}(u)) \right] 
  \ - \ 
  \lambda_{s,*} \left( \left[ \w_\fg(u) , \w_\fg(v) \right] \right) \, .
\end{eqnarray}
In particular, $(\lambda_{p+1,1} , \lambda_{p,2} , \ldots , \lambda_{p+2-k,k})$ must lie in the kernal of the map
\begin{displaymath}
  \delta_{p+2} : \bigoplus \fg^\perp_{p+2-m} \ot \fg_{-m}^* \ \to \ 
  \bigoplus  \fg^\perp_{p+2-i-j} \ot ( \fg_{-i}^* \wedge \fg_{-j}^* )
\end{displaymath}
defined as follows:  given $x \in \fg_{-i}$ and $y \in \fg_{-j}$,
\begin{equation}\label{eqn:sp2la}
\renewcommand{\arraystretch}{1.3}
\begin{array}{l}
  \delta_{p+2} (\lambda_{p+1,1} , \lambda_{p,2} , \ldots , \lambda_{p+2-k,k})(x,y)
  \ =     \\
  \hbox{\hspace{100pt}}
  \left[ x , \lambda_{p+2-j,j}(y) \right] \ - \ 
  \left[ y , \lambda_{p+2-i,i}(x) \right] \ - \ 
  \lambda_{p+2-i-j,i+j} \left( \left[ x,y \right] \right) \, .
\end{array}\end{equation}
As we will see below, $\delta_{p+2} = \partial^1_{p+2}$ is the degree p+2 component of the first Lie algebra cohomology differential (described in \S\ref{liealgcohsect}).  As in \S\ref{sec:warmup}, $\delta_{p+2}$ should be viewed as the Spencer differential augmented to deal simultaneously with
tableau and  torsion.


Moreover, $\fg^\perp_{p+2}  = \fn\cap  \fgl(U)_{p+2}$ (cf. \S\ref{redprolongsect}), and the Lie algebra cohomology denominator $\partial^0_{p+2}(\fg^\perp_{p+2})$ is the space of admissible normalizations of the prolongation coefficients $\lambda_{p+2-m,m}$.  Thus, the vanishing of $H^1_{p+2}(\fg_-,\fg^\perp)$ implies that normalized integral manifolds of the $(I^\textsf{f}_{p},\Omega)$ system are in one to one
correspondence with integral manifolds of the $(I^\textsf{f}_{p+1},\Omega)$ system.
We have established the following lemma.

\begin{lemma} \label{translemma}
Let $U$ be a complex vector space, and $\fg\subset \fgl(U)$ a represented complex semi-simple Lie algebra.  Let $ G/P \subset \bP U$ be the corresponding homogeneous variety.  Denote the induced $\bZ$-gradings by $\fg = \fg_{-k} \op \cdots \op \fg_k$ and $U = U_0 \op \cdots \op U_{-f}$.  Fix an integer $p \ge -1$, and let $(I_p^\textsf{f},\Omega)$  denote the
filtered    linear Pfaffian system  given by
 \eqref{psystem}.

If $H^1_{p+2}(\fg_{-},\fg\upperp)=0$,   then  normalized integral manifolds
of the $I^\textsf{f}_p$ system are in one to one correspondence with integral manifolds
of the $I^\textsf{f}_{p+1}$ system.
\end{lemma}

\noindent Now a straightforward induction yields our main result.

\begin{theorem} \label{transthm}\label{mainthm}
Let $U$ be a complex vector space, and $\fg\subset \fgl(U)$ a represented complex semi-simple Lie algebra.  Let $Z = G/P \subset \bP U$ be the corresponding homogeneous variety  (the orbit of a highest weight line).  Denote the induced $\bZ$-gradings by $\fg = \fg_{-k} \op \cdots \op \fg_k$ and $U = U_0 \op \cdots \op U_{-f}$.  Fix an integer $p \ge -1$, and let $(I_p^\textsf{f},\Omega)$  denote the linear Pfaffian system given by \eqref{psystem}.
If $H^1_d(\fg_{-},\fg\upperp)=0$, for all $d\geq p+2$, then the homogenous variety $G/P$ is rigid for the $(I_p^\textsf{f},\Omega)$ system.
\end{theorem}

\section{The $(I_p^\textsf{f},\Omega)$-systems and Lie algebra cohomology}\label{cohrulessect}
\subsection{Lie algebra cohomology}\label{liealgcohdef}
\label{liealgcohsect}
Let $\fl$ be a Lie algebra and let $\Gamma$ be an $\fl$-module.
The maps
$$
\partial^{j}: \La j\fl^*\ot \Gamma \ra \La{j+1}\fl^*\ot \Gamma
$$
are defined in a natural way to respect the Leibnitz rule \cite{kostant}, and give rise to a complex.
Define $H^k(\fl,\Gamma):=\tker\partial^k/\tim \partial^{k-1}$.
We will only have need of $\partial^0$ and $\partial^1$ which are   as follows.  If $X\in \Gamma$ and $v,w\in \fl$, then
$$\partial^0(X)(v)=v.X, 
$$
 and 
  if $\a\ot X\in \La 1\fl^*\ot \Gamma$,  then
\begin{equation}\label{del1}
\partial^1(\a\ot X)(v\ww w)=
\a([v,w])X+\a(v)w.X - \a(w) v.X \, .
\end{equation}

Now let  $\fl$ be a graded Lie algebra and $\Gamma$   a graded
$\fl$-module. The
chain complex and Lie
algebra cohomology groups inherit  gradings as well.
Explicitly,
$$
\partial^{1}_{d}: \oplus_{i}(\fl_{-i})^* \ot \Gamma_{d-i}
\ra \oplus_{j\leq m}(\fl_{-j})^*\ww (\fl_{-m})^*\ot \Gamma_{d-j-m} \, .
$$
Taking
$\fl=\fg_{-}$, $\Gamma=\fg\upperp$ and $d=p+2$, we obtain $\d_{p+2}=\partial^1_{p+2}$ 
in the $I^\textsf{f}_p$ system as asserted following (\ref{eqn:sp2la}).  Moreover, the image of $\partial^0_{p+2}$ is the space of admissible normalizations.

\subsection{Applying Kostant's theory}\label{kosapply}\label{kossect}
Kostant \cite{kostant} shows that under the following circumstances one can
compute $H^k(\fl,\Gamma)$ combinatorially:
\begin{enumerate}
\item $\fl=\fn\subset\fp\subset  \fg$ is the nilpotent subalgebra of
a parabolic subalgebra of
a semi-simple Lie algebra $\fg$.

\item $\Gamma$ is a $\fg$-module.
\end{enumerate}

Under these conditions, letting $\fg_0\subset \fp$ be the the
(reductive) Levi factor of $\fp$, $H^k(\fn,\Gamma)$ is naturally
a $\fg_0$-module.  When $\Gamma$ is an irreducible $\fg$-module of
highest weight $\l$, the irreducible $\fg_0$-modules appearing in $H^k(\fn,\Gamma)$ have highest weight $w.(\lambda)$, $w \in \cW^\fp(k)$. Here, $\cW^\fp(k)$ is a subset of the Weyl group of $\fg$, and $w.(\l)=w(\l+\rho)-\rho$ denotes the affine action of $w$, with $\rho=\sum_i\o_i=\frac 12\sum_i\a_i$ the sum of the fundamental weights of $\fg$.

For the $ I_p^\textsf{f}$ system associated to $ G/P_I$ we will only be concerned with $H^1(\fg_{-},\fg\upperp)$.  In this case $\cW^\fp(1)$ is the set of simple reflections $\sigma_i$ corresponding to simple roots $\a_i$ such that $i \in I$.    We will use the grading element $Z$, introduced in \S\ref{rootgradingsect}, to decompose the first cohomology group into homogeneous components $H^1(\fg_- , \fg^\perp) = \bigoplus H^1_d(\fg_-,\fg^\perp)$. 
\medskip

\noindent{\it Remark.}     Our $\fg_-$ is Kostant's $\fn^* = \fl^*$.  In particular, $H^1_d(\fg_-,\Gamma) = H^1_{-d}(\fn,\Gamma)$; and it is the latter that we shall be computing via Kostant.

\subsection{Computing \boldmath $H^1_{d}(\fg_-,\fg\upperp)$ \unboldmath}\label{compsect}

By \S\ref{kosapply}, we need only calculate the affine
reflections in the simple roots $\a_{i}$
and to determine if they are nonzero
in degree greater than $p+1$. We follow \cite{HY,SYY}
in this subsection.

Let $\sigma_{i_0}$ denote the affine reflection through simple root $\a_{i_0}$, $i_0 \in I$.  If we apply $\s_{{i_{0}}}$
to $\l=\l^i\o_i$ we obtain, from Equation \eqref{charelteqn},
$$
Z(\s_{\a_{i_{0}}}.\l) \ = \ \sum_{i,s}\l^i(c\inv)_{i,i_{s}}-\l^{i_{0}}-1.
$$

To see this, let 
$\langle \cdot,\cdot\rangle: \ft^*\times \ft^*\ra \BC$
be the pairing determined by $\langle\o_i,\a_j\rangle=\d_{ij}$.
\begin{eqnarray*}
Z(\s_{\a_{i_0}}.\l)&= &
Z(\s_{\a_{i_0}}(\l+\rho))-Z(\rho)
\\
&=&Z(\l+\rho -\langle \l+\rho,\a_{i_0}\rangle\a_{i_0})-
Z(\rho)
\\
&=&Z(\l)+Z(\rho) -\langle \l+\rho,\a_{i_0}\rangle Z(\a_{i_0})
-Z(\rho)
\\
&=&Z(\l) -\langle \l+\rho,\a_{i_0}\rangle 
\\
&=& \sum_{i,s} \l^i(c\inv)_{i,i_s}- (\l_{i_0}+1)
\end{eqnarray*}
In particular, for the trivial representation,
$Z(\s.0)=-1$.

For positive weights $\l$, $Z(\s_{\a_{i_0}}.\l) > -1$ when $(c\inv)_{i_{0},i_{0}}>1$; and $Z(\s_{\a_{i_0}}.\l) \ge -1$ when $(c\inv)_{i_{0},i_{ 0}}\ge1$.
Recall, our $H^1_d(\fg_-,\Gamma)$ is Kostant's $H^1_{-d}(\fn,\Gamma)$ (cf. \S\ref{kossect}).  Thus, 
\smallskip
\begin{center}
{\it $H^1_d(\fg_-,\Gamma) = 0$, for $d\ge1$ (resp. $d>1$) when $(c^{-1})_{i_0,i_0} > 1$ (resp. $\geq 1$), for all $i_0 \in I$.}
\smallskip
\end{center}

The inverse Cartan matrices satisfy $(c\inv)_{j,j}>1$,
except for the following cases (omitting redundancies):
\begin{eqnarray*}
 (c\inv)_{1,1} \, = \, (c^{-1})_{n,n} \ = \ \frac n{n+1}  
 &{\rm \ for\ }&  A_n \\
 (c\inv)_{2,2} \ = \ 1  &{\rm \ for\ }&  A_3 \\
 (c\inv)_{1,1} \ = \ 1  &{\rm \ for\ }&  B_n, \ n\geq 2, \ \textrm{ and } \ D_n, \ n\geq 4 \, .
\end{eqnarray*}
By Schur's lemma, if   $\fg$ is a semi-simple Lie algebra
and $V$ is an irreducible $\fg$-module, then
the trivial representation occurs exactly
once in $\fgl(V)=V^*\ot V$ and not at all in $\fsl(V)$.
These observations, coupled with Theorem \ref{transthm}, allow us to deduce Theorem \ref{bigthm} and 

\begin{theorem}\label{thm:0rigidity} 
Let $G$ be a complex semi-simple group and let $G/P\subset \BP U$ be a homogeneously embedded homogeneous variety. Assume that $G/P$ contains no factor corresponding a homogeneous variety $A_n/P_I$, with $1$ or $n$ in $I$.  Then, the $(I_0^\mathsf{f},\Omega)$ system on $SL(U)$ is rigid. 
\end{theorem}
\medskip

\noindent{\it Proof of Theorem \ref{adjointthm}.}
By Proposition \ref{prop:Fub3=>-1} integral manifolds for the third-order Fubini system
for adjoint varieties  are also integral manifolds of the $(I_{-1},J_{-1})$ system and hence the $(I_{-1}^f,\Omega)$ system. 
Theorem \ref{adjointthm} now follows from Theorem \ref{bigthm} and Lemma \ref{adjointlemma} below. 
\hfill\qed
\smallskip
\begin{lemma}\label{adjointlemma}
The $(I_{-1}^f,\Omega)$ system
is rigid for   $U=U^{A_n}_{ \o_1+\o_n }$, $n>2$.
 The $(I_0^\textsf{f},\Omega)$ system
is rigid for $U=U^{A_2}_{ \o_1+\o_2 }$.
\end{lemma}

\noindent{\it Remark.}  This establishes Theorem \ref{anadjthm}.

\begin{proof}
We have the decomposition,
$$U^*\ot U=U_{\o_1+\o_n}\ot U_{\o_1+\o_n}=
U_{2(\o_1+\o_n)}\op U_{2\o_1+\o_{n-1}}
\op U_{\o_2+2\o_n}
\op U_{\o_2+\o_{n-1}}
\op 2U_{\o_1+\o_n}
\op U_0
$$
so
$$\fg\upperp=
U_{2(\o_1+\o_n)}\op U_{2\o_1+\o_{n-1}}
\op U_{\o_2+2\o_n}
\op U_{\o_2+\o_{n-1}}
\op  U_{\o_1+\o_n}   \, .
$$

We calculate, assuming $n\geq 2$ (note that when $n=2$, $\o_{n-1}=\o_1$)
\begin{eqnarray*}
Z(\s_{\a_{i_0}}.2(\o_1+\o_n))&=& 1, \quad i_0=1,n \\
Z(\s_{\a_1}.(2 \o_1+\o_{n-1}))&=&  \left\{\begin{matrix}
1 & n\geq 3\\  -1 & n=2\end{matrix} \right. \\
Z(\s_{\a_1}.(  \o_2+2\o_{n }))&=&  \left\{\begin{matrix}
2 & n\geq 3\\  2 & n=2\end{matrix} \right. \\
Z(\s_{\a_1}.(  \o_2+ \o_{n-1 }))&=&  \left\{\begin{matrix}
1 & n\geq 3\\  0 & n=2\end{matrix} \right. \\
Z(\s_{\a_{i_0}}. (\o_1+\o_n))&=& 0, \quad i_0=1,n \, .
\end{eqnarray*}
The $ Z(\s_{\a_n}.\l)$'s may be determined
by symmetry.
\end{proof}


Now we consider some examples of the $(I_0^\textsf{f},\Omega)$ system.

\medskip

\noindent
{\bf Example:} {\it Fubini's theorem.}
The third order Fubini system for a quadric hypersurface is equivalent
to the $(I_0^\textsf{f},\Omega)$ system once one observes
$\tann (\mathit{II})=\fg_0+\BC\subset \fgl(U)_0$.  We obtain a calculation free proof 
of Fubini's theorem.
\smallskip

Rigidity of the $(I_0^\textsf{f},\Omega)$ system can only fail for $A_n/P_I$ if $\fg^\perp$ contains a module with large $\w_1$ (resp. $\w_n$) coefficient and $1 \in I$ (resp. $n\in I$), and  all other $\w_j$ coefficients are relatively small.  Theorem \ref{vdpna2thm} is an example for which rigidity holds.  
\medskip

\noindent{\bf Proof of Theorem \ref{vdpna2thm}.}
For $v_d(\pp n)\subset\BP S^d\BC^{n+1}=\BP U$ we consider
$U\ot U^*=U_{d\o_1}\ot U_{d\o_n}=\sum_{i=0}^dU_{i(\o_1+\o_n)}$.
Hence $\fg\upperp=\sum_{i=2}^dU_{i(\o_1+\o_n)}$.
We calculate $Z(\s_{\a_1} .i(\o_1+\o_n))=i-1$.

A similar argument establishes the second half of the theorem for $U^{A_2}_{\w_1+\w_2}$.  (The rigidity of $U^{A_n}_{\w_1+\w_n}$, $n>2$, is given by Lemma \ref{adjointlemma}.)
\hfill\qed\medskip

We show in \S\ref{sl3sect} that the nontrivial
Lie algebra cohomology group $H^1_1$ that occurs
from  $Z(\s_{\a_1}.(3\o_1))$ (here $\o_{n-1}=\o_1$)
indeed makes the $(I_{-1} , J_{-1})$ system flexible.
 
\section{Flexibility of ${\tSeg}(\pp 1\times \pp n)$}\label{p1pnsect}

Here we prove Theorem \ref{thm:flexseg}.  Fix vector spaces $E$ and $F$ of dimensions $2$ and $n+1$, respectively.  Let $\tSeg( \bP E \times \bP F) \subset \bP(E\ot F)$ denote the Segre variety.  If $n=1$, then $\tSeg(\bP^1\times\bP^1)$ is a quadric hypersurface; quadrics
are  flexible at order two and rigid at order three by Fubini's theorem \cite{fub}.  So we assume for the rest of this section that $n>1$. We have already seen above
that ${\tSeg}(\pp 1\times \pp n)$ is rigid to order three
because it is rigid for the $(I_0,J_0)$ system.

We consider the $(I_{-1},J_{-1})$ system which agrees the second order Fubini system  in this case.  Here
$W=\gp{-1}$,
$A=\gp 0$, and
$V=\fg_{-1}$.

However,
 $H^1_1(\fg_{-},\fg\upperp)\simeq \BC^{n}$, which
as a $\fg_0=\BC\op (\fsl_{n}+\BC)$ module
is acted on by weights $-4$ for the first $\BC\subset\fsl_2$
and $1$ for the second $\BC\subset \fa_n$
and is the dual of the standard representation for
the  $\fsl_{n}$ factor, corresponding
to the weight
$-4\eta + \o^1+\o^n$ where $\eta$ is the fundamental
weight of $\fa_1$ and $\o_j$ are fundamental weights
of $\fa_n$.  

We now take our EDS to be the reduced prolongation.
At this point the Spencer differential still coincides
with the Lie algebra cohomology group, but the
tableau no longer corresponds to a $\fg$-module,
so we do not see any obvious way to use Kostant's
theory. Fortunately in this case the calculation
is straightforward and the prolonged system is involutive,
with characters $(s_1,s_2 \hd s_{n+1})=(n^2+2n,0\hd 0)$.  Moreover, it is not difficult to check that any $Y \subset \bP^{2n+1}$ agreeing with $\tSeg$ to second order admits a (reduced) frame bundle $\cF \subset \tGL_{2n+2}\bC$ on which the Maurer-Cartan form pulls-back to
$$
\begin{pmatrix}
\ooo 00 &\ooo 0 1&\ooo 0\b&r^0_{\ob}\oo 1\\
\oo 1 &\ooo 11& 0 &\ooo 0\b\\
\oo\a & r^{\a}_{1}\, \oo 1&\ooo\a\b& r^{\a}_{\ob}\oo 1
+\d^{\a}_{\b}\ooo 01\\
0&\oo\a&\d^{\a}_{\b}\oo 1&\ooo\a\b+\d^{\a}_{\b}(\ooo 11-\ooo 00)
\end{pmatrix} \, .
$$
Here $2 \le \a \le n+1$, and $\oa = \a + n$, and forms $\w^0_0$, $\w^1_1$, 
$\w^\a_\b$, $\w^1_0$, $\w^\a_0$, $\w^0_1$ and $\w^0_\b$ are linearly independent.  The $r^\a_{1}$ terms were introduced in the first reduced
prolongation, and the $r^0_\ob$ and $r^\a_\ob$ coefficients in the second. 

Given an integral manifold of the system, let $Y \subset \bP^{2n+1}$ denote the corresponding base manifold.  By construction, $Y$ agrees with $\tSeg(\bP E \times \bP F)$ to second order.  Third order agreement holds if and only if the Fubini cubic vanishes on $Y$.  This is equivalent to $r^\a_1 = 0$.  This completes
the proof of Theorem \ref{thm:flexseg}. 
\medskip
 
\noindent{\it Remark.}
All (base) integral manifolds $Y$ of our system are ruled by $\bP^n$s.  Consider the induced curve in the Grassmannian 
$\mathbb{G}(\bP^n,\bP^{2n+1})$. Such curves agree with the curve induced by $\tSeg(\bP E \times \bP F)$ to order two, but differ with this curve
at order three.  In a remark in the introduction in \cite{Lchss} it was mistakenly stated that the above system was rigid, the error was primarily caused by the misconception that a curve in the Grassmannian would be determined by its second fundamental form.

\section{The $(I_{-1},J_{-1})$ system for $A_2$}\label{sl3sect}
Here we consider $Y^3 \subset \bP^7$ with frame bundle $\cF_{<0}$ on which the 
$\mathfrak{sl}(\fa_2)_{<0}$ component of the Maurer-Cartan form $\w$ agrees 
with that of the $\mathfrak{sl}(3,\bC)$ adjoint variety $X$.  Explicitly,  
\begin{equation} \label{eqn:sl3neg} \renewcommand{\arraystretch}{1.3}
  \w \ = \ \left( 
  \begin{array}{cccccccc}
    \w^0_0 & \b_1 & \b_2 & \b_3 & \w^0_4 & \w^0_5 & 
    \w^0_6 & \w^0_7 \\
    \cline{1-1}
    \multicolumn{1}{c|}{\a^1} & \w^1_1 & \w^1_2 & \w^1_3 & \w^1_4 & 
    \w^1_5 & \w^1_6 & \w^1_7 \\
    \multicolumn{1}{c|}{\a^2} & \w^2_1 & \w^2_2 & \w^2_3 & \w^2_4 & 
    \w^2_5 & \w^2_6 & \w^2_7 \\ \cline{2-3}
    \a^3 & -\half \a^2 & 
    \multicolumn{1}{c|}{\half \a^1} & \w^3_3 & \w^3_4 & \w^3_5 & \w^3_6 
    & \w^3_7 \\ 
    0 & -\a^2 & 
    \multicolumn{1}{c|}{-\a^1} & \w^4_3 & \w^4_4 & \w^4_5 & \w^4_6 & 
    \w^4_7 \\ 
    \cline{4-5}
    0 & -\a^3 & 0 & -\a^1 & \multicolumn{1}{c|}{\threehalves \a^1} & 
    \w^5_5 & \w^5_6 & \w^5_7 \\
   0 & 0 & -\a^3 & -\a^2 & \multicolumn{1}{c|}{-\threehalves \a^2} & \w^6_5 &
    \w^6_6 & \w^6_7 \\ \cline{6-7}
   0 & 0 & 0 & -2\a^3 & 0 & -\a^2 & \multicolumn{1}{c|}{\a^1} & \w^7_7
  \end{array} \right) \, .
\end{equation}
The $\w^0_a = \b_a$, are not assignments, merely renamings.

The   quadrics of the second fundamental form are 
\begin{eqnarray} \label{eqn:sl3quadrics}
  q^4 \ = \ -2 \, \a_1 \, \a_2 & & q^5 \ = \ -2 \, \a_1 \, \a_3 \\
  q^6 \ = \ -2 \, \a_2 \, \a_3 & & q^7 \ =\  -2 \, \a_3{}^2 \, . \nonumber
\end{eqnarray}
The tableau is torsion-free with characters $(s_1,s_2,s_3) = (17,4,0)$.  In the course of prolonging, the 1-forms $\w^1_2$ and $\w^2_1$ are forced to be semi-basic (linear combinations of the $\a^a$).  There exist rescalings of the $e_j$, $e = (e_0 , \ldots , e_7 ) \in \cF$, that preserve the EDS and scale the $\a^2$--coefficient of $\w^1_2$ and the $\a^1$--coefficient of $\w^2_1$ to be constant.  We decompose the analysis into three cases: that the constants are $\{0,0\}$, $\{0,1\}$ (two symmetric cases) or $\{1,1\}$.

\subsection{Case 1} \label{sec:case1}
Assume that both constants are zero.    Taking the second prolongation shows that $\w$ satisfies the $(I_0,J_0)$ system. Thus,  integral manifolds of (\ref{eqn:sl3neg}) are frame bundles over the adjoint variety $Z^{A_2}_\tad$.
\subsection{Case 2}\label{sec:case2}  
Suppose that $\w^1_2 \equiv \a^2$ mod $\a^1, \a^3$; and $\w^2_1 \equiv 0$ mod $\a^2, \a^3$.  After restricting our parameter space (twice) to remove torsion, the second tableau has characters $(s_1 , s_2 , s_3) = (13,1,0)$.  The (reduced) second prolongation is of dimension 9.  Again, we must restrict the parameter space (twice, again) to obtain a torsion-free third tableau with characters $(s_1,s_2,s_3) = (9,1,0)$.  The third prolongation has dimension 5.  The fourth tableau (after one reduction for torsion) has characters $(s_1,s_2,s_3) = (4,0,0)$ and is  involutive.  The Maurer-Cartan form (the components of non-negative degree) is given below.
In degree 0 (the block diagonal):
\footnotesize
\begin{eqnarray*}
  \w^0_0 & = & f_1 \, \a^1 + f_2 \, \a^2 + \fourth f_{1,2} \, \a^3 + 3 \, \w^1_1 
  \, , \quad  \w^1_2  \ = \ \a^2 \, , \quad 
  \w^5_6 \ = \  \a^2 - \half \, f_1 \, \a^3 \\ 
  0 & = & \w^0_0 - \w^1_1 - \w^2_2 \ = \ \w^2_2 + \w^5_5 \ = \  
  \w^1_1 + \w^6_6 \ = \ \w^0_0 + \w^7_7 \ = \ \w^2_1 \ = \ \w^3_3 \ = \ 
  \w^3_4 \ = \ \w^4_3 \ = \ \w^4_4 \ = \ \w^6_5 \, .
\end{eqnarray*}
\normalsize
Above the $f_j$ are functions on the (reduced) frame bundle, and the 
$f_{j,a}$ are the $\a^a$--coefficients of $\td f_j$.  In degree 1:
\footnotesize
\begin{displaymath}\renewcommand{\arraystretch}{1.3}
\begin{array}{r@{ \ \ = \ \ }l@{\qquad}r@{ \ \ = \ \ }l}
  \w^1_3 & \w^4_6 \ = \ \fourth \, f_1 \a^2 - \b_2 \, , &
  \w^1_4 & \tfrac98 \, f_1 \, \a^2 - \tfrac32 \, \b_2 \, , \quad
  \w^3_6 \ = \ -\tfrac18 \, f_1 \, \a^2 - \fourth \, f_{1,1} \, \a^3 
    - \half \, \b_2 \\
  \w^5_7 & -\fourth \, f_{1,1} \, \a^3 - \b_2 \, , &
  \b_1 & \w^2_3 \ = \ -\tfrac23 \, \w^2_4 \ = \ -2 \, \w^3_5 \ = \ 
  \w^4_5 \ = \ \w^6_7 \, .
\end{array}\end{displaymath}
\normalsize
In degree 2: 
\footnotesize
\begin{displaymath}\renewcommand{\arraystretch}{1.3}
\begin{array}{r@{ \ = \ }l@{\qquad}r@{ \ = \ }l@{\qquad}r@{ \ = \ }l}
  \w^0_4 & -\tfrac38 \, f_{1,1} \, \a^2 - \tfrac32 \, p_1 \, \a^3 &
  \w^1_5 & \tfrac18 \, f_{1,1} \, \a^2 - \tfrac32 \, p_1 \, \a^3 
    - \half \, \b_3 &
  \w^1_6 & \fourth \, f_{1,2} \, \a^2 + \fourth \, f_{1,3} \, \a^3 - \b_1 \\
  
  \w^4_7 & p_1 \, \a^3 \, , \quad  \w^2_5 \, = \, 0 \, &
  \w^2_6 & -\tfrac18 \, f_{1,1} \, \a^2 - \tfrac32 \, p_1 \, \a^3 
    - \half \, \b_3 &
  \w^3_7 & -\tfrac18 \, f_{1,1} \, \a^2 -6 \, p_1\, \a^3 - \half \, \b_3 \, , 
\end{array}\end{displaymath}
\normalsize
with $p_1 = (f_1 \, f_{1,1} + f_{1,1,1})/24$, where $f_{1,1,1}$ is the 
$\a^1$--coefficient of $\td f_{1,1}$.  Note that $p_1$ appears below in the Fubini cubic $r^4$.
In degree 3: 
\footnotesize
\begin{displaymath}\renewcommand{\arraystretch}{1.3}
\begin{array}{r@{ \ = \ }l@{\qquad}r@{ \ = \ }l}
  \w^0_5 & -4 \, p_1 \, \a^2 - p_3 \, \a^3 & 
  \w^0_6 & -p_1 \, \a^1 - \half \, f_{1,3} \, \a^2 - 2 \, p_2 \, \a^3 
    - \tfrac34 \, f_1 \, \b_1 \\
  \w^1_7 & -p_1 \, \a^1 + \fourth \, f_{1,3} \, \a^2 + p_2 \, \a^3 
    + \fourth \, f_1 \, \b_1 &
  \w^2_7 & -p_1 \, \a^2 + p_3 \, \a^3 \, , 
\end{array}\end{displaymath}
\normalsize
with $p_2 = (f_{1,1} \, f_{1,2} + 4 \, f_{1,1,3} )/48$ and 
$p_3 = ( f_{1,1}{}^2 - 2 \, f_1{}^2 \, f_{1,1} - 2 \, f_1 \, f_{1,1,1} )/48$.
And the final degree 4 1-form is 
\footnotesize
\begin{displaymath}
  \w^0_7 \ = \ 
  -4 \, p_3 \, \a^1 - 4 \, p_2 \, \a^2 + p_4 \, \a^3 
  - \half \, f_{1,1} \,\b_1 \, , 
\end{displaymath}
\normalsize
where $p_4$ is a polynomial in $f_1 , f_2 , f_{1,a} , f_{1,1,1}, f_{1,1,3}$ (homogeneous of degree 4 in derivatives).

The Fubini cubics are
\begin{eqnarray*}
  r^4  \ = \ -2 \, \a_2{}^3 + 2 \, p_1 \, \a_3{}^3 & & 
  r^5 \ = \ 3 \, \a_1{}^2 \, \a_2 - \tfrac32 \, f_1 \, \a_2 \, \a_3{}^2 
    - \half \, f_{1,1} \, \a_3{}^3 \\
  r^6 \ = \ -3 \, \a_1 \, \a_2{}^2 & & r^7 \ = \ 0 \, .
\end{eqnarray*}

These varieties are intrinsically flat (cf. \S\ref{sec:contact}).

\subsection{Case 3} \label{sec:case3}
Assume that $\w^1_2 = \a^2$ mod $\a^1, \a^3$; and $\w^2_1 = \a^1$ mod $\a^2, \a^3$.
The second tableau (after two sets of restrictions on our parameter space to remove torsion) has characters $(s_1,s_2,s_3)=(15,2,0)$.  The second (reduced) prolongation has dimension 14.  The third tableau (after one set of restrictions to remove torsion) has characters 
$(s_1,s_2,s_3) = (11,2,0)$.  The third prolongation is (also) of dimension 14.  The fourth tableau (no torsion) has characters $(s_1,s_2,s_3)=(12,2,0)$, and the fourth prolongation is of dimension 16.  Thus the system is involutive by Cartan's test.

The Maurer-Cartan form (the components of non-negative degree) is given below.
In degree 0:
\footnotesize
\begin{eqnarray*}
  \w^1_1 & = & f_1 \, \a^1 + f_2 \, \a^2 + p_1 \, \a^3 \, , \quad 
  \w^2_2 \ = \ g_1 \, \a^1 + g_2 \, \a^2 + p_2 \, \a^3 \, , \quad
  \w^1_2 \ = \ \a^2 \, , \quad \w^2_1 \ = \ \a^1 \\
  \w^5_6 & = & \a^2 + \tfrac12 \, ( 2 \, f_1 - g_1 ) \, \a^3 \, , \quad 
  \w^6_5 \ = \ \a^1 + \tfrac12 \, ( f_2 - 2 \, g_1 ) \, \a^3 \\
  0 & = & \w^0_0 - \w^1_1 - \w^2_2 \ = \ \w^1_1 + \w^6_6 \ = \ 
  \w^2_2 + \w^5_5 \ = \ \w^3_3 \ = \ \w^3_4 \ = \ \w^4_3 \ = \ \w^4_4 \ = \ 
  \w^0_0 +  \w^7_7  \, .
\end{eqnarray*}
\normalsize
The $f_j$, $g_j$, $1 \le j \le 2$, are functions on the frame-bundle; the $f_{j,a}$ and $g_{j,a}$ are the $\a^a$--coefficients of $\td f_j$ and $\td g_j$; and
$\td f_1 \equiv \b_1$, $\td f_2 \ \equiv \ -2 \, \b_2$, 
$\td g_1 \ \equiv \ -2 \, \b_1$ and $\td g_2 \ \equiv \ \b_2$ 
mod the semi-basic $\a^a$.  Finally $p_1 = -1 + f_1 \, ( f_2 - g_2 ) + f_{2,1} - f_{1,2}$ and $p_2 = 1 + ( f_1 - g_1 ) \, g_2 + g_{2,1} - g_{1,2}$. 
In degree 1:
\footnotesize
\begin{displaymath}\renewcommand{\arraystretch}{1.3}
\begin{array}{r@{ \ = \ }l@{\qquad}r@{ \ = \ }l@{\qquad}r@{ \ = \ }l}
  \w^1_3 & -\tfrac14 \, ( 2 \, f_1 - g_1 ) \, \a^2 - \b_2 & 
  \w^1_4 & -\tfrac98 \, ( 2 \, f_1 - g_1 ) \, \a^2 - \tfrac32 \, \b_2 &
  \w^3_5 & -\tfrac18 \, ( f_2 - 2 \, g_2 ) \, \a^1 + p_3 \, \a^3 
    - \half \, \b_1 \\
  \w^2_4 & \tfrac98 \, ( f_2 - 2 \, g_2 ) \, \a^1 - \tfrac32 \, \b_1 &
  \w^2_3 & -\tfrac14 \, ( f_2 - 2 \, g_2 ) \, \a^1 + \b_1 & 
  \w^3_6 & \tfrac18 \, ( 2 \, f_1 - g_1 ) \, \a^2 - p_4 \, \a^3 
    - \half \, \b_2 \\
  \w^4_5 & -\tfrac14 \, ( f_2 - 2 \, g_2 ) \, \a^1 + \b_1 & 
  \w^4_6 & -\tfrac14 \, ( 2 \, f_1 - g_1 ) \, \a^2 - \b_2 &
  \w^5_7 & p_4 \, \a^3 - \b_2 \quad
  \w^6_7 \, = \, -p_3 \, \a^3 + \b_1 \, .
\end{array}\end{displaymath}
\normalsize
Above, $p_3 = \tfrac 14 \, ( 2 \, f_1 - g_1 ) + \tfrac12 \, g_2 \, ( f_2 - 2 \, g_2 ) + \tfrac14 \, ( f_{2,2} - 2 \, g_{2,2} )$ and $p_4 = \tfrac14 \, ( f_2 - 2 \, g_2 ) + \tfrac12 \, f_1 \, ( 2 \, f_1 - g_1 ) + \tfrac14 \, ( 2 \, f_{1,1} - g_{1,1} )$.
In degree 2:
\footnotesize
\begin{displaymath}
\begin{array}{c}
   \w^0_4 =  -\tfrac32 \, p_3 \, \a^1 + \tfrac32 \, p_4 \, \a^2 
     + h_1  \, \a^3 \, ,  \quad 
   \w^3_7 =  \half \, p_3 \, \a^1 + \half \, p_4 \, \a^2 
     + h_2 \, \a^3 - \half \, \b_3 \, , \quad
   \w^4_7 =  -\tfrac23 \, h_1 \, \a^3 \\
   \w^1_5 = \tfrac14 \, (2 \, f_1 - g_1 + 2 \, p_3 ) \, \a^1 
     - \tfrac14 \, ( f_2 - 2 \, g_2 + 2 \, p_4 ) \, \a^2 + h_1 \, \a^3
     - \half \, \b_3  \\
   \w^1_6 = \tfrac14 \, (f_2 - 2 \, g_2 ) \, \a^1 - ( 2 \, p_1 + p_2 ) \, \a^2 
     + p_5 \, \a^3 - \b_1 \, , \quad
   \w^2_5 = (p_1 - 2 \, p_2) \, \a^1 + \tfrac14 \, ( 2 \, f_1 - g_1 ) \, \a^2 
     + p_6 \, \a^3 + \b_2 \\
   \w^2_6 = -\tfrac14 \, ( 2 \, f_1 - g_1 + 2 \, p_3 ) \, \a^1 
     + \tfrac14 \, ( 2 \, f_2 - g_2 + 2 \, p_4 ) \, \a^2 + h_1 \, \a^3
     - \half \, \b_3 \, .
\end{array}
\end{displaymath}
\normalsize
In degree 3:
\footnotesize
\begin{displaymath}
\begin{array}{rcl}
  \w^0_5 & = & 2 \, p_6 \, \a^1 + \left( \tfrac83 \, h_1 
    + \tfrac{3}{16} \, (2 \, f_1 - g_1) (f_2 - 2 \, g_2) \right) \, \a^2 
  + p_7 \, \a^3 - \tfrac34 \, ( f_2 - 2 \, g_2 ) \, \b_2 \\
  \w^0_6 & = & \left( \tfrac23 \, h_1 
    - \tfrac{3}{16} \, (2 \, f_1 - g_1) (f_2 - 2 \, g_2) \right) \, \a^1 
    - 2 \, p_5 \, \a^2 + p_8 \, \a^3  + \tfrac34 ( 2 \, f_1 - g_1 ) \, \b_1 \\
  \w^1_7 & = & \left( \tfrac23 \, h_1 
    + \tfrac{1}{16} \, (2 \, f_1 - g_1) (f_2 - 2 \, g_2) \right) \, \a^1 
    + p_5 \, \a^2 + p_9 \, \a^3 - \tfrac14 \, ( 2 \, f_1 - g_1 ) \, \b_1 \\
  \w^2_7 & = &  p_6 \, \a^1 + \left( \tfrac43 \, h_1 
    + \tfrac{1}{16} \, (2 \, f_1 - g_1) (f_2 - 2 \, g_2) \right) \, \a^2
    + p_{10} \, \a^3 + \tfrac14 \, (f_2 - 2 \, g_2) \, \b_2
\end{array}
\end{displaymath}
\normalsize
The $p_7 , \ldots , p_{10}$ above are polynomials in $\{ f_j , g_j , h_j , f_{j,k} , g_{j,k} , h_{j,k} \}_{1\le j,k\le 2}$.
Finally, in degree 4:
\footnotesize
\begin{displaymath}
  \w^0_7 \ = \ p_{11} \, \a^1 + p_{12} \, \a^2
    + \left( \half \, h_{2,3} + (p_1 + p_2 ) \, h_2 \right) \, \a^3
    - \tfrac14 \, ( f_2 - 2 \, g_2 - 8 \, p_4 ) \, \b_1
    + \tfrac14 \, ( 2 \, f_1 - g_1 + 4 \, p_3 ) \, \b_2 \, .
\end{displaymath}
\normalsize
The $p_{11} , p_{12}$ above are polynomials in the $\{ f_j , g_j , h_j , f_{j,k} , g_{j,k} , h_{j,k} \}_{1\le j,k\le 2}$.  

Finally, we remark that the functions $f_j , g_j , h_j$ satisfy a system of partial differential equations.  For example, $( f_2 + 2 \, g_2 )_1 = ( 2 \, f_1 + g_1 )_2$.  One should think of this system as reducing the freedom of our initial data from four functions of $14$ variables
to the two functions of two variables and   twelve functions of one variable and  specified by the characters.

The Fubini cubics are
\begin{eqnarray*}
  r^4  & = & -2 \, ( \a_1{}^3 + \a_2{}^3 ) 
    - \tfrac43 \, h_1 \, \a_3{}^3 \\
  r^5 & = &  3 \, \a_1{}^2 \, \a_2 
    + \tfrac32 \, ( 2\, f_1 - g_1 ) \, \a_2 \, \a_3{}^2 
    + 2 \, p_4 \, \a_3{}^3 \\
  r^6 & = & -3 \, \a_1 \, \a_2{}^2 
    + \tfrac32 \, ( f_2 - 2 \, g_2 ) \, \a_1 \, \a_3{}^2 
    - 2 \, p_3 \, \a_3{}^3 \\ 
  r^7 & = &  0 \, .
\end{eqnarray*}
In $r^4$, $h_1$ is the term that first appeared in $\w^0_4$ above.

The intrinsic invariants (cf. \S\ref{sec:contact}) are 
\begin{eqnarray*}
  J_1 & = & - \tfrac18 \, ( 2 \, f_1 - g_1 ) - \tfrac12 \, ( 3 \, f_2 - g_2 ) \, g_2
    - \tfrac14 \, f_{2,2} + \tfrac43 \, p_3 \\
  J_2 & = & -\tfrac18 \, ( f_2 - 2 \, g_2 )  - \tfrac12 \, f_1 \, ( f_1 + 2 \, g_1 )
    - \tfrac12 \, f_{1,1} - \tfrac13 \, p_4 \, .
\end{eqnarray*}
\subsection{Intrinsic invariants}
\label{sec:contact}
Let $Y^3 \subset \bP^{7}$ be a 3-fold with frame-bundle $\cF_\tneg$ on which the 
Maurer-Cartan form is as given in (\ref{eqn:sl3neg}).  From (\ref{eqn:sl3quadrics}) we see that $\BP
\{ v\in T_xX\mid F_2(v,v)=0\} \subset \bP T_{-1,x}$ consists of exactly two points.  These two 
points define a pair of line bundles which span the contact hyperplane $T_{-1,x} \subset T_xX$.

Let's consider the situation more generally.  Let $M$ be a complex 3-dimensional manifold (not necessarily projective) admitting two line bundles 
$L_1$ and $L_2$ such that $H := L_1 \op L_2$ is a contact distribution.
Consider the bundle $\cF \to M$ of all frames $\{ e_1 , e_2 , e_3 \}$ of 
$T_xM$ such that $e_j$ spans $L_j$, and given any local section 
$e : U \subset M \to \cF$, $[e_1(x) , e_2(x)] \equiv e_3(x)$ mod $H_x$.
Then $\cF$ is a principle $G$-bundle with fibre group $G \subset \tGL_3\bC$ 
\begin{displaymath}
  G \ := \ \left\{ \, \left.
    \left( \begin{array}{ccc}
     a & 0 & \a \\
     0 & b & \b \\
     0 & 0 & ab 
   \end{array} \right) \, 
   \right| \, 
     \begin{array}{l}
     a , \, b \, \in \bC \, , \\ 
     \a , \,  \b \, \in \bC\backslash\{0\} 
     \end{array} \, \right\} \, .
\end{displaymath}

Given $e = \{ e_j \} \in \cF_x$ and $v\in T_e\cF$, the canonical semi-basic 
$\bC^3$-valued 1-form on $\cF$ is given by $\pi_*(v) = \eta^j(v) e_j$.
There exist connection 1-forms $\{ \theta^s \}_{s=1}^4$ on $\cF$ such that 
\begin{eqnarray*}
  \td \eta^1 & = & \eta^1 \wedge \th^1 + \eta^3 \wedge \th^3 \\
  \td \eta^2 & = & \eta^2 \wedge \th^2 + \eta^3 \wedge \th^4 \\
  \td \eta^3 & = & -\eta^1 \wedge \eta^2 + \eta^3 \wedge ( \th^1 + \th^2 ) \, .
\end{eqnarray*}
Notice that this connection has torsion (see $\td \eta^3$).  In fact, there 
are no torsion-free connections on $\cF$.  The family of connections 
preserving the structure equations above is 4-dimensional.

The choice of connection may be further refined (leaving one degree of 
freedom) so that 
\begin{eqnarray*}
  \td \th^1 & = & 2 \, \th^4 \wedge \eta^1  + \th^3 \wedge \eta^2 
                  - \phi \wedge \eta^3 \\
  \td \th^2 & = & - \th^4 \wedge \eta^1 - 2 \, \th^3 \wedge \eta^2 
                  - \phi \wedge \eta^3 \\
  \td \th^3 & = & \th^2 \wedge \th^3 - \phi \wedge \eta^1 
                  + J_1 \eta^2 \wedge \eta^3 \\
  \td \th^4 & = & \th^1 \wedge \th^4 - \phi \wedge \eta^2 
                  + J_2 \eta^1 \wedge \eta^3 \, .
\end{eqnarray*}
Above, $\phi$ is a 1-form, and the $J$ are functions on $\cF$.

The functions $J$
are relative invariants. They define invariant tensors
   $J_1 \, (\w^2 \wedge \w^3) \ot \ule^3, J_2 \, (\w^1 \wedge \w^3) \ot \ule^3\in \Gamma(M, \bigwedge^2 (T^*M) \ot (TM/H)^*)$.

Applying the above
formulas to Case 2 of \ref{sec:case2}
shows the two invariants $J_1$ and $J_2$ are zero, so the integral manifolds are {\it intrinsically} flat. 
Similarly, we see  Case 3 (\S\ref{sec:case3}) is not intrinsically flat.

\subsection{Order three rigidity of \boldmath $Z^{A_2}_{\tad}$\unboldmath}\label{sl3end}
 It is also true  that $ I_0\subset I_{\tFub_3}$ but
this is only seen to hold after one calculates
a derivative. One can either prove third
order rigidity of $X_\tad^{A_2}$ this way
or cite \cite{robles}.

\section{Remaining rigidity proofs}\label{restsect}

\subsection{The varieties \boldmath$v_d(\pp n)$, $v_d(Q^n)$ and
$\tSeg(\pp{1}\times \pp{a_2}\ctimes \pp{a_n})$\unboldmath}\label{resta}  Assume throughout this section that $n>1$.

\medskip

\noindent{\it Proof of Theorem \ref{thm:I0J0} (a).}  Begin with the Veronese variety, and consider the Fubini system of order $d+1$. We will
show it is {\it a priori} more restrictive than the $(I_0,J_0)$ system;
thus the $(I_0^\textsf{f},\Omega)$ system being rigid implies
it is rigid.
Here $N_k\simeq S^kT^*$, the $F_{k,k}$ are the
identity maps for $k\leq d$ and all the $F_{s,k}=0$ for $s<k$
and $F_q=0$ for $q>d$.
Since the grading is three step, the osculating filtration
coincides with the Lie algebra filtration (albeit with
different integers attached to the filtrands).

The vanishing of the $F_{k,k-1}$, $k=3\hd d$ fixes
$N^*_k\ot N_k$ in terms of $T^*\ot T$ and $L^*\ot L$,
thus together they fix
$\fgl(U)_{0}$ to be $\fg_0$.
The vanishing of the $F_{k,k-2}$, $k=4\hd d$ fixes
$N^*_k\ot N_{k-1}$ in terms of 
$L^*\ot T$ and the remaining forms in the system
fix all other components of the Maurer-Cartan
form taking values in the remaining spaces
below the diagonal to be zero.
Thus we are reduced to the $(I_0,J_0)$ system, and Theorem \ref{thm:I0J0} (a) follows from Theorem \ref{vdpna2thm}. \hfill \qed
\medskip

\noindent{\it Remark.} Note that we only used a small part of the Fubini
system to prove rigidity in this case. 
\medskip 

\noindent{\it Remark.} In in the case $d=2$, $n=1$, $v_2(\bP^1) \subset \bP^2$ is a plane conic.  Monge showed that the plane conics are rigid at order 5.  This is consistent with the Lie algebra cohomology:  In this case we have $\fg = U = S^2\bC^2 = U_{2 \w_1}$.  The decomposition $\fsl(U) = U_{2\w_1} \op U_{4 \w_1}$ yields $\fg^\perp = U_{4 \w_1}$.  The cohomology group $H^1(\fg_- , \fg^\perp)$ is one-dimensional, with weight $-6\w_1$ as a $\fg_0 = \bC$ module.  Since the grading element $Z$ for $U$ is given by $Z(\w_1) = \half$, we see that $H^1_3(\fg_-,\fg^\perp) \not=0$, an obstruction to Fubini rigidity at order four.
\medskip

\noindent{\it Proof of Theorem \ref{thm:I0J0} (b).}
The $v_d(Q)$ Fubini system is similar
to the Veronese system. All Fubini forms
are zero except for the fundamental forms
and they are given by  
$N_k\simeq S^kT^*$ for $k\leq d$, 
$N_{d+1}=Q\circ S^{d-1}T^*$, 
$N_{d+2}=Q^2\circ S^{d-2}T^*$,..., $N_{2d}=Q^d$ and all higher fundamental
forms are zero.
Here $F_{2d+1,2d}$ fixes the $N_d^*\ot N_d$ component
of the Maurer-Cartan form in terms of the
$T^*\ot T\op L^*\ot L$-component, and
similarly down the line. Thus 
Theorem \ref{thm:I0J0} (b) follows from Theorem \ref{bigthm}.
\hfill \qed

\medskip

\noindent{\it Proof of Theorem \ref{thm:I0J0} (c).}  We merely sketch the proof for the Segre variety $\tSeg(\pp 1\times \pp{a_2}\ctimes\pp{a_r})$ -- it is similar to the two above.  The $r+1$ Fubini system implies the $(I_0,J_0)$ system.  We compute  that
$H^1_d(\fg_{-},\fg\upperp)=0$ for $d>1$.  Thus the $(I_0,J_0)$ system
is rigid. \hfill\qed

\subsection{\boldmath$E_8/P_1\subset \pp{3874}$\unboldmath}\label{e8sect}

To establish   that $E_8/P_1$ is rigid to order five,
it suffices to show that the osculating sequence has length four
and that the $\o_{\gp{< 0}}$ component
of the Maurer-Cartan form is forced
to vanish on integral manifolds
of the fifth order Fubini system.
The first assertion follows because the shortest sequence of negative roots taking $\w_1$ to $-\w_1$ is of length four.
The second by arguments similar to the adjoint case
for the $T^*_{-1}\ot T_{-1}$ component, forcing
this component of the Maurer-Cartan form to be the
spin representation of $D_{7}$ and
by observing that the fourth fundamental form is
the quadric on $T_{-2}$, forcing
the $T^*_{-2}\ot T_{-2}$ component of the
Maurer-Cartan form to be the standard representation of $D_7$.



\end{document}